\documentclass[11pt,reqno]{conm-p-l}

\usepackage{amsmath}
\usepackage{amsfonts}
\usepackage{amscd}
\usepackage{enumerate}
\usepackage{fancyheadings}
\usepackage{bbold}

\usepackage{amsmath}
\usepackage{amsthm}
\usepackage{amssymb}
\usepackage{amsfonts}
\usepackage{epsfig}
\usepackage{epsf}
\newtheorem{theorem}{Theorem}[section]

\newtheorem{EVtheorem}{Etingof - Varchenko Theorem} %in ch4
\newtheorem*{Stheorem}{Schiffmann Theorem} %in ch4

\newtheorem*{remark}{Remark}
\newtheorem*{acknow}{Acknowledgments}
\newtheorem{lemma}{Lemma}[section]

\usepackage{latexsym}
\usepackage{epsfig}

\newcommand{\baseRing}[1]{\ensuremath{\mathbb{#1}}}
\newcommand{\C}{\baseRing{C}}
\newcommand{\Z}{\baseRing{Z}}

\renewcommand{\P}{\baseRing{P}}

\newcommand{\g}{\mathfrak{g}}
\newcommand{\h}{\mathfrak{h}}
\renewcommand{\b}{\mathfrak{b}}

\newcommand{\Spec}{\operatorname{Spec}}

\renewcommand{\phi}{\varphi}
\allowdisplaybreaks
\begin{document}

\title{Dynamical Quantum Groups - The Super Story}
\author{Gizem Karaali}
\address{Department of Mathematics, Pomona College, Claremont, CA
91711}
\email{gizem.karaali@pomona.edu}

\begin{abstract}
We review recent results in the study of quantum groups in the super setting.
In particular, we provide an overview of results about solutions of the
Yang-Baxter equations in the super setting and develop the super analog of the
theory of dynamical quantum groups.
\end{abstract}

\maketitle

\section{Introduction}
\label{SectionIntroduction}

\subsection{Background and Goals}
\label{BackgroundandGoals}

Algebraists have been studying Hopf algebras since the 1960's. Written first in 1969, \cite{Swe} is still one of the classical references in the subject. The quantum Yang-Baxter equation (QYBE) was already known in 1950s, at least within the mathematical physics community. \cite{Ji2} contains some of the earliest papers on the subject. However, the discovery of the particular noncommutative noncocommutative Hopf algebras called \emph{quantum groups} relating these two concepts has undoubtedly intensified the research in both areas over the last two decades. Along the way, a full theory of quantum groups has been developed, to the extent that there are currently various textbooks on the subject. See, for instance, \cite{BG, CP, ES, HK, Kas, Maj}.

Super structures have been of interest to mathematicians and physicists alike. Besides providing the mathematical framework for supersymmetry, they have
proved to be mathematically rich structures.

In this paper we collect together results and ideas that can be helpful in the
pursuit of a full theory of quantum groups in the super setting. Currently,
there are only partial results in this direction.
Although we do mention some new results, our main goal is to provide a clear
exposition of the present state of affairs in the theory of super quantum
groups with a strong emphasis on dynamical structures. In particular, we focus on the  solutions of the classical and quantum dynamical Yang-Baxter equations, with the standard theory for the non-dynamical equations as our guide.

For Sections \ref{SectionDynamicalrmatrices} and \ref{SectionSuperDynamicalrmatrices}, we will assume some familiarity with the theory of Lie superalgebras. \cite{Kac} has the first comprehensive study of these structures. \cite{Sche} and \cite{Var} provide some relevant background on the subject.\footnote{Another beautiful reference with a more geometric flavor is \cite{Ma97}. Here, Manin develops projective algebraic geometry in the context of supermanifolds. In this note, we will need neither the strength nor the sophistication of the tools provided there.} A concise summary of sign conventions used in the study of super structures can be found in \cite{DeFr}. For Sections \ref{SectionDynamicalQuantumGroups} and \ref{SectionSuperGroups}, some familiarity with quantum groups and Hopf algebras at the level of a text like \cite{Kas} will be sufficient.

\begin{acknow}
The author thanks P. Etingof, L. Feh{\'{e}}r, A. Isaev, E. Koelink, M.
Kotchetov, J-H. Lu, S. Montgomery, H. Rosengren and R. Wisbauer for
suggestions and constructive comments during the work that led to this paper, and the two referees whose recommendations improved this paper significantly. 
It is also a pleasure to thank the organizers L. Kauffman, D. Radford and
F. Souza, of the AMS Special Session on Hopf Algebras at the Crossroads of
Algebra, Category Theory, and Topology, October 23-24, 2004,
where the author had the opportunity to present her results in \cite{Kar2}.
\end{acknow}

\subsection{Plan of this Paper}
\label{PlanofPaper}

Section \ref{SectionIntroduction} is introductory.. In \S\ref{SuperSolutionsNonDynamical}, we give a brief overview of the
results from \cite{Kar1} and \cite{Kar2}. We explain their relevance to our
ultimate goal without going into too much detail. In 
\S\ref{MotivationforDynamicalfirst}, we explain our motivation for the emphasis
of this paper on the dynamical picture.

In the first (classical) part of the paper, consisting of Sections
\ref{SectionDynamicalrmatrices} and \ref{SectionSuperDynamicalrmatrices}, we
study the classical dynamical Yang-Baxter equation (CDYBE) and its
solutions. In \S\ref{HistoricalOverview} we give a brief
overview of the historical development of the subject of dynamical Yang-Baxter
equations. \S\ref{DefinitionsDynamical} provides the precise
definitions of the terms involved. \S\ref{SchiffmannResult} is a
concise but explicit summary of the classification results for the non-graded
case. Section \ref{SectionSuperDynamicalrmatrices} is concerned with various
super analogues for the results from Section \ref{SectionDynamicalrmatrices}.
We make the appropriate definitions in
\S\ref{DefinitionsSuperDynamical}, and present some construction and
classification results for the super solutions of the CDYBE in
\S\ref{SuperZeroWeight} and \S\ref{GeneralizingSchiffmann}. The material in
\S\ref{SuperZeroWeight} appeared elsewhere \cite{Kar3}, but the main result
of \S\ref{GeneralizingSchiffmann} is new.

Sections \ref{SectionDynamicalQuantumGroups} and \ref{SectionSuperGroups} make
up the second (quantum) part of the paper. In Section
\ref{SectionDynamicalQuantumGroups}, we describe the general theory of
dynamical quantum groups. In order to provide a comprehensible exposition, we
begin, in \S\ref{WhatisHopfAlgebroid}, with the definitions of
groupoids, bialgebroids, and Hopf algebroids. In \S\ref{QuantumGroupoids}, we
discuss dynamical quantum groups in more detail. We consider the categorical picture in 
\S\ref{CategoricalSubsection}. In Section \ref{SectionSuperGroups}, we begin our study of the super analogue of the theory of dynamical quantum groups. We consider the super versions of the basic definitions in \S\ref{DefinitionsSuperGroups}. In
\S\ref{ResultsSuperGroups}, we discuss the QDYBE and its solutions in the
super setting.

Section \ref{SectionConclusion} closes the present exposition with a brief
discussion of some open problems and a possible plan
of action for them.

\subsection{Super Solutions of the Classical Yang-Baxter Equation}
\label{SuperSolutionsNonDynamical}

The classical Yang-Baxter equation (CYBE) for an element $r\in{\g}\otimes \g$
where $\g$ is a Lie algebra with a nondegenerate $\g$-invariant bilinear form
is: 
\[ [r^{12}, r^{13}]+[r^{12}, r^{23}]+[r^{13}, r^{23}]=0.\footnote{Here, and elsewhere in the paper, the notation $r^{12}$ stands for the sum $\sum a_i \otimes b_i \otimes 1$ where $r = \sum a_i \otimes b_i$. Similarly, $r^{13} = \sum a_i \otimes 1 \otimes b_i$ and $r^{23} = \sum 1 \otimes a_i \otimes b_i$.} \]
\noindent
In \cite{BD1} and \cite{BD2}, Belavin and Drinfeld classified the solutions of
the CYBE, which are called $r$-matrices. In particular, their work provides us
with an explicit construction of non-skewsymmetric $r$-matrices associated to
certain discrete data (the Belavin-Drinfeld triples).

In \cite{Kar1}, the author proved a similar construction result in the case of
simple Lie superalgebras $\g$ with nondegenerate forms. Hence, the
constructive part of the result of Belavin and Drinfeld mentioned above was
extended to the framework of superalgebras. However, the discrete data that
directly generalize the Belavin-Drinfeld triples were seen to be insufficient
to fully classify all (non-skewsymmetric) solutions to the CYBE in the super
setting. In particular, the author introduced (in \cite {Kar1}) and studied in
detail (in \cite{Kar2}) a particular non-skewsymmetric $r$-matrix on $sl(2,1)$
that cannot be distinguished from the standard $r$-matrix on the basis of
the Belavin-Drinfeld type data alone. Such examples, in fact, can be
constructed in any $sl(m,n)$ for any $m \not = n$.

As is well-known, the solutions to the CYBE on a Lie (super)algebra give us the
semiclassical limits of quantum $R$-matrices (solutions to the QYBE) on the associated Lie (super)group.
Hence, an understanding of the solutions of the CYBE is essential to the theory
of super quantum groups that we would like to develop. So far the results
mentioned above show clearly that there will be some intrinsically new
structures that need to be considered, although perhaps a mere change in
perspective will suffice to see these as natural extensions of the original
non-graded results.\footnotemark

\footnotetext{This comment is due to I. Penkov.}

\subsection{Why Study Dynamical Structures?}
\label{MotivationforDynamicalfirst}

Quantum $R$-matrices provide us with examples of quantum groups.
The standard quantum
$R$-matrix is the quantization of the standard $r$-matrix
obtained from the trivial Belavin-Drinfeld triple. The natural problem of
quantizing all of the Belavin-Drinfeld type $r$-matrices was solved completely
first by Etingof and Kazhdan in \cite{EK}. However, their method was very abstract, so the explicit constructions by Etingof, Schedler and
Schiffmann in \cite{ESS} were a most welcome development. This latter method
in fact works for the dynamical $r$-matrices, i.e. the solutions of the more
general classical dynamical Yang-Baxter equation (CDYBE). The associated
quantum objects in this dynamical setting are the dynamical quantum groups (or
dynamical quantum groupoids, to be more precise), and the associated Hopf
objects are in fact Hopf algebroids.

In this paper, we concentrate on the super solutions of the dynamical
Yang-Baxter equation and the super versions of the notions of dynamical
quantum groups and Hopf algebroids. Naturally we expect that understanding
these will help us in extending the quantization
result from \cite{ESS} cited above to obtain a graded analogue. 
It should be noted here that the method of Etingof and Kazhdan has already
been generalized to the super setting, in \cite{Ge}. This is a very important
development. However, due to the less constructive and more abstract nature of
\cite{EK} and \cite{Ge}, we think that following the dynamical
route can still provide us with additional valuable insight.

\section{Dynamical \textit{r}-matrices and the Classical Dynamical Yang-Baxter
Equation}
\label{SectionDynamicalrmatrices}

\subsection{Historical Overview}
\label{HistoricalOverview}

Even before the classical dynamical Yang-Baxter
equation (CDYBE) first appeared, the quantum dynamical Yang-Baxter equation (QDYBE) had been
considered by Gervais and Neveu in \cite{GN}. Their motivation
had been purely physical. Indeed, they were studying monodromy matrices
in Liouville theory. 

In order to talk about the introduction of the CDYBE, we need to briefly
mention conformal field theory. (For more information on the following material, we refer the reader to \cite{EFK, Kas}). The correlation functions in
Wess-Zumino-Witten (WZW) conformal field theory can be constructed out of
holomorphic sections of certain vector bundles; these sections are known as
\emph{conformal blocks}. The 
%correlation functions of 
conformal blocks on ${\P}^1$ for
the WZW conformal field theory for a simple Lie
algebra $\g$ satisfy the Knizhnik-Zamolodchikov (KZ) equations,
which were introduced in \cite{KZ}, and generalized in
\cite{Ch}. KZ equations play an important role in representation
theory and mathematical physics. For a system of KZ equations to
have a solution, the necessary and sufficient consistency
condition coincides with the classical Yang-Baxter equation
(\cite{Ch, EFK}). 

The conformal blocks for an arbitrary elliptic curve satisfy a system of
differential equations, a modification of the KZ equations,
known as the Knizhnik-Zamolodchikov-Bernard (KZB)
equations (\cite{Bern}). The classical Yang-Baxter equation
is no longer the required consistency condition for the KZB
equations. However, it turns out that a modification will work.

Felder introduced the CDYBE formally in \cite{Fe1} and \cite{Fe2}. In fact, it
had been appearing in the mathematical physics literature in disguise for a
while. (See \cite{BDF}, for instance,  where some of its basic trigonometric
solutions were studied in detail). In \cite{Fe1} and \cite{Fe2}, Felder showed
that the CDYBE was the consistency condition for the KZB
equations. He also related it, in a way analogous to the non-dynamical case,
to the quantum dynamical Yang-Baxter equation (QDYBE). 

At this point, it may be worth mentioning that, in this context, the label \emph{dynamical} refers to the fact that the relevant
equation is now a differential rather than an algebraic
equation, and so it may remind us of a dynamical
system.

The classical and quantum dynamical Yang-Baxter equations, as well as their
various solutions, have been extensively studied. In the late 1990s, Etingof and
Varchenko
started to develop the classification of solutions of these
equations, and together with Schiffmann they completed a
classification theory for the classical case under certain
assumptions (\cite{ES, EV, Schi}). 
Their classification results resemble the classification results of
Belavin and Drinfeld (\cite{BD1, BD2}) in the
non-dynamical case.

As Drinfeld pointed out in \cite{D}, solutions of the
classical Yang-Baxter equation have a natural geometrical
interpretation as Poisson-Lie structures on the
corresponding Lie groups. Similarly,
solutions of the CDYBE
define Poisson-Lie groupoid structures.
See \cite{EV, FM, KS1, KS2, Lu2, Wei} for research in this direction.

By now, the theory of the classical and quantum
dynamical Yang-Baxter equations and their solutions has many
applications, including some to integrable systems and
representation theory. For a recent survey of results and such
applications, one can refer to \cite{Et, ES2}. The text \cite{EtLa} provides a very readable introduction to the topic as well.

\subsection{Definitions}
\label{DefinitionsDynamical}
Let $\g$ be a simple Lie algebra and ${\h \subset \g}$ a
Cartan subalgebra. The \emph{classical dynamical Yang-Baxter
equation} for a meromorphic function ${r : \h^* \rightarrow \g
\otimes \g}$ is 
\begin{equation}
\label{DYBequation}
Alt(dr) + 
[r^{12},r^{13}] + [r^{12},r^{23}] + [r^{13},r^{23}] = 0;
\end{equation}
\noindent
where the differential of $r$ itself is considered as a
meromorphic function:  
$$ \begin{matrix}
dr &:& \h^* &\longrightarrow& \g \otimes \g \otimes \g  \\
&& \lambda &\longmapsto& \sum_i x_i \otimes \frac{\partial
r}{\partial x_i} (\lambda)  .
\end{matrix} $$
\noindent
Here $\{x_i\}$ is a basis for $\h$, but it is also regarded as a
linear system of coordinates on $\h^*$. 

Using this we have:
$$ Alt(dr) = 
\sum_i x_i^{(1)} \left(\frac{\partial r}{\partial x_i}
\right)^{(23)} +  
\sum_i x_i^{(2)} \left(\frac{\partial r}{\partial x_i}
\right)^{(31)} +  
\sum_i x_i^{(3)} \left(\frac{\partial r}{\partial x_i}
\right)^{(12)}.$$ 
%alt is the symmetrization wrt even permutations of 1,2,3
\noindent
Therefore, we can rewrite the classical dynamical Yang-Baxter
equation as follows: 
$$\sum_i x_i^{(1)} \left(\frac{\partial
r}{\partial x_i} \right)^{(23)} 
+ \sum_i x_i^{(2)} \left(\frac{\partial r}{\partial x_i}
\right)^{(31)} 
+ \sum_i x_i^{(3)} \left(\frac{\partial r}{\partial x_i}
\right)^{(12)} + [[r,r]] = 0. $$
%$$+ [r^{12},r^{13}] + [r^{12},r^{23}] + [r^{13},r^{23}] = 0 
%$$

A meromorphic function $r : \h^* \rightarrow \g \otimes \g$ is
called a \emph{dynamical r-matrix associated to $\h$} if it is a
solution to the classical dynamical Yang-Baxter equation
\eqref{DYBequation} and it satisfies both the \emph{zero weight condition} and the \emph{generalized unitarity condition}. These two conditions are described as follows:

\begin{enumerate}
\item The \emph{zero weight condition} for a meromorphic function $r : \h^* \rightarrow \g \otimes \g$ is: 
$$ [h\otimes 1 + 1 \otimes h, r(\lambda)] = 0 \textmd{ for all }
h \in \h, \lambda \in \h^*.$$
\item \emph{Generalized unitarity}\footnotemark 
\footnotetext{This is a generalization of the \emph{unitarity condition} which
corresponds to $\epsilon = 0$. If $R$ is a quantization of $r$, then the
statement $R T(R) = 1$ can be considered as a quantization of the statement $r
+ T(r) = 0$, and hence the term \emph{unitarity}. At this point, it may be useful to note that, in this paper, we do not focus on the actual process of quantization, but we always keep in mind the relationship between the classical and the quantum pictures. This allows us to make such interpretations as in this footnote, while still keeping this paper at a reasonable length.}
with coupling constant
$\epsilon$ for a meromorphic function $r : \h^* \rightarrow \g \otimes \g$ is:  
\begin{equation}
\label{dynamicalunitarity} 
r(\lambda) + T(r)(\lambda) = \epsilon \Omega.
\end{equation}
\end{enumerate}
\noindent
Here $\Omega$ is the Casimir element, i.e. the element of ${\g
\otimes \g}$ corresponding to the Killing form. $T$ is the standard
\emph{twist} map on the second tensor power of any vector space: $T(a\otimes
b) = b \otimes a$.

\subsection{Classification Results}
\label{SchiffmannResult}

Below is a brief survey of recent results related to the
classification of the solutions of the classical dynamical
Yang-Baxter equation.

The first classification results about the solutions of the classical dynamical Yang-Baxter equation can be found in \cite{EV}. There, Etingof and Varchenko proved two theorems that result in a full classification of all dynamical $r$-matrices satisfying the zero weight condition:
%\vspace{0.1in}
\begin{EVtheorem}
\label{EVtheorem1}
(1) Let $X$ be a subset of the set of roots $\Delta$ of a
simple Lie algebra $\g$ with nondegenerate Killing form ${(
\cdot \; , \cdot )}$ such that:

(a) If $\alpha, \beta \in X$ and $\alpha + \beta$ is a root,
then $\alpha + \beta \in X$, and

(b) If $\alpha \in X$, then $-\alpha \in X.$

\noindent
Let $\nu \in \h^*$, and let $D = \sum_{i<j} D_{ij} dx_i \wedge
dx_j$ be a closed meromorphic $2-$form on $\h^*$. If we set
$D_{ij} = -D_{ji}$ for $i \ge j$, then the meromorphic
function:   
$$ r(\lambda) =
\sum_{i,j=1}^N D_{ij}(\lambda) x_i \otimes x_j +  \sum_{\alpha
\in X} \frac{1}{(\alpha,
\lambda-\nu)} e_{\alpha} \otimes e_{-\alpha} $$ 
\noindent
is a dynamical r-matrix with zero weight and zero coupling
constant.

\noindent
(2) Any dynamical r-matrix with zero weight and zero
coupling constant is of this form. 
\end{EVtheorem}

\begin{EVtheorem} 
\label{EVtheorem3}
(1) Let $\g$ be a simple Lie algebra with
nondegenerate Killing form ${( \cdot \; , \cdot )}$. Let $\Delta$ be 
the set of roots of $\g$, and fix a subset $X$ of the set $\Gamma$
of simple positive roots. Let $\overline{X}$ be the intersection
of the linear span of $X$ with $\Delta$. Pick a $\nu \in
\h^*$, and define: 
$$\phi_{\alpha} = \left \{ \begin{matrix}
\left(\epsilon / 2\right)\coth\left(
\left(\epsilon / 2\right) (\alpha,\lambda - \nu) \right) 
& \textmd{   if  } \alpha \in \overline{X} \\
\left(\epsilon / 2\right) 
& \textmd{   if  } \alpha \not \in \overline{X}, 
\textmd{  positive} \\  
-\left(\epsilon / 2\right) 
& \textmd{   if  } \alpha \not \in \overline{X},
\textmd{  negative}  
\end{matrix} \right. $$
\noindent
Let $D = \sum_{i<j} D_{ij} dx_i \wedge dx_j$ be a closed
meromorphic $2-$form on $\h^*$. If we set $D_{ij} = -D_{ji}$ for
$i \ge j$, then the meromorphic function: 
$$ r(\lambda) =
\sum_{i,j=1}^N D_{ij}(\lambda) x_i \otimes x_j +
\frac{\epsilon}{2}\Omega + \sum_{\alpha \in \Delta}
\phi_{\alpha} e_{\alpha} \otimes e_{-\alpha}  $$   
\noindent
is a dynamical r-matrix with zero weight and nonzero
coupling constant $\epsilon$. 

\noindent
(2) Any classical dynamical r-matrix with zero weight
and nonzero coupling constant $\epsilon$ is of this form. 
\end{EVtheorem}

These are Theorems $3.2$, and $3.10$, respectively in \cite{EV}. 
We will provide super analogs of these theorems in \S\ref{SuperZeroWeight}.

A more general theory for classifying dynamical r-matrices, with no need for the zero weight condition, has been developed by Schiffmann in \cite{Schi}. Below is a brief summary of his results. Later in \S\ref{GeneralizingSchiffmann}, we will discuss their super versions.

The general setup of \cite{Schi} is as follows:
Let $\g$ be a simple Lie algebra with nondegenerate invariant
bilinear form ${(\cdot,\cdot)}$, ${\mathfrak{l} \subset \g}$ a
commutative subalgebra containing a regular semisimple element
on which ${(\cdot,\cdot)}$ is nondegenerate \footnotemark, $\h$ the Cartan
subalgebra containing $\mathfrak{l}$, and $\h_0$ the orthogonal
complement of $\mathfrak{l}$ in $\h$. Denote by $\Omega_{00}$
the $\h_0-$part of the Casimir element $\Omega$ of $\g$. Fix a set $\Gamma$ of simple roots. 
\footnotetext{This condition about the existence of a regular element was later
shown to be redundant (\cite{Et}).}

In this setup, Schiffmann starts by introducing the notion of a
a \emph{generalized Belavin-Drinfeld triple}, which is defined as a triple
${(\Gamma_1, \Gamma_2, \tau)}$ where both ${\Gamma_i}$
are subsets of the set $\Gamma$ of simple roots, and ${\tau :
\Gamma_1 \rightarrow \Gamma_2}$ is a norm-preserving bijection.
Then, an
$\mathfrak{l}-$\emph{graded} generalized Belavin-Drinfeld triple
is one where $\tau$ preserves the decomposition of $\g$ into
$\mathfrak{l}-$weight spaces. Given a generalized
Belavin-Drinfeld triple ${(\Gamma_1, \Gamma_2, \tau)}$, denote
by $\Gamma_3$ the largest subset of $\Gamma_1 \cap \Gamma_2$
that is stable under $\tau$. Clearly $\tau$ can be extended to
an isomorphism ${\g_{\Gamma_1} \rightarrow \g_{\Gamma_2}}$ once
we fix a set of Weyl-Chevalley generators.

If we fix basis vectors ${\{e_{\alpha}\}}$ for the
non-Cartan part of $\g$, as is done in the standard Belavin-Drinfeld
construction, then for any ${\lambda \in \mathfrak{l}^*}$, we
can define a map ${K(\lambda) : \mathfrak{n}_+(\Gamma_1)
\rightarrow \mathfrak{n}_+(\Gamma_2)}$ as follows:
\begin{eqnarray*} 
K(\lambda) (e_{\alpha}) &=& \frac{1}{2} e_{\alpha} +
e^{-(\alpha,\lambda)}
\frac{\tau}{1-e^{-(\alpha,\lambda)}\tau}(e_{\alpha}) \\
&=& \frac{1}{2}e_{\alpha} + 
\sum_{n>0} e^{-n(\alpha,\lambda)}\tau^n(e_{\alpha}). 
\end{eqnarray*}
\noindent
Note that the sum above is finite as long as ${\alpha \not\in
\Gamma_3}$. If ${\alpha \in \Gamma_3}$ and $\tau|_{\Gamma_3} = id_{\Gamma_3}$,
then we have:
$$ K(\lambda) (e_{\alpha}) = \frac{1}{2} 
\coth\left(\frac{1}{2}(\alpha,\lambda) \right)e_{\alpha}.$$

Given a meromorphic map ${r : \h^*  \rightarrow (\g
\otimes \g)^{\h}}$, we consider the following transformations:
\begin{enumerate}
\item $r(\lambda) \longmapsto r(\lambda) + \sum_{i<j}
C_{ij}(\lambda) x_i \wedge x_j$, where $\sum_{i,j}C_{ij}(\lambda)
d \lambda_i \wedge d \lambda_j$ is a closed meromorphic
$2-$form,  
\item $r(\lambda) \longmapsto r(\lambda -\nu)$, where $\nu \in
\h^*$, 
\item $r(\lambda) \longmapsto (A\otimes A) \; r(A^*\lambda)$,
where $A$ is an element of the Weyl group of $\g$.
\end{enumerate}
\noindent
We will call these transformations \emph{gauge transformations}. 
Two dynamical r-matrices that can be obtained from one another
via a sequence of gauge transformations are called
\emph{gauge-equivalent}.

The main result of \cite{Schi} is the following:

\vspace{0.1in}
\begin{Stheorem} 
(1) Any dynamical r-matrix will be gauge-equivalent to a dynamical
r-matrix $r$ satisfying: 
\begin{equation}
\label{lskewsymmetry}
r(\lambda) - T(r)(\lambda) \in \mathfrak{l}^{\perp} \otimes 
\mathfrak{l}^{\perp} = \left( (\bigoplus_{\alpha \neq 0}
\g_{\alpha} ) \oplus \h_0 \right) \otimes
\left((\bigoplus_{\alpha \neq 0} \g_{\alpha}) \oplus
\h_0\right)
\end{equation}
\noindent
Here, once again, $T$ is the usual twist, 
mapping $a\otimes b$ to $b \otimes a$.

\noindent
(2) Let ${(\Gamma_1, \Gamma_2, \tau)}$ be an
$\mathfrak{l}-$graded generalized Belavin-Drinfeld triple. 
If ${r_{00} \in \h_0 \otimes \h_0}$ satisfies: 
\begin{equation*}
(\tau(\alpha) \otimes 1) r_{00}  + (1 \otimes \alpha)r_{00} =
\frac{1}{2} ((\alpha + \tau(\alpha)) \otimes 1) \Omega_{00}
\end{equation*}
\noindent
for any $\alpha \in \Gamma_1$, then:
$$ r(\lambda) = \frac{1}{2}\Omega + r_{00} + \sum_{\alpha \in
\overline{\Gamma}_1 \cap \Delta^+} K(\lambda)(e_{\alpha}) \wedge
e_{-\alpha} + \sum_{\alpha \in \Delta^+, \alpha \not \in
\overline{\Gamma}_1} \frac{1}{2}e_{\alpha} \wedge e_{-\alpha}$$
\noindent
is a dynamical r-matrix satisfying Equation
\eqref{lskewsymmetry}.

\noindent
(3) Any dynamical r-matrix satisfying Equation
\eqref{lskewsymmetry} is of the above form for a suitable
triangular decomposition of $\g$. 
\end{Stheorem}

\vspace{0.1in}

\begin{remark}
Note that when $\mathfrak{l}= \h$ we get Etingof - Varchenko
Theorem \ref{EVtheorem3}. The two subsets $\Gamma_1, \Gamma_2$
are the subset $X$ of the mentioned theorem, and the two
statements are equal when we set $\epsilon = 1$.  The
isometry $\tau$ is the identity on $X$. The zero-weight
condition is equivalent to the associated generalized
Belavin-Drinfeld triple being $\mathfrak{l}-$graded.
\end{remark}

\vspace{0.1in}

\begin{remark}
If $\Gamma_3 = \emptyset$, then ${(\Gamma_1, \Gamma_2, \tau)}$
is a Belavin-Drinfeld triple. If $\tau$ preserves the grading of
the $\mathfrak{l}-$weight space decomposition of $\g$, then 
$\mathfrak{l}$ has to be orthogonal to the
coroots corresponding to the roots in the span of $(\Gamma_1
\cup \Gamma_2)$. If ${\mathfrak{l} = \{0\}}$, then
${\mathfrak{l}^* = \{0\}}$, and we get constant
(non-dynamical) r-matrices, i.e. solutions of the CYBE. The statement in this
case
can be seen to be equivalent to the Belavin-Drinfeld
Theorem.
\end{remark}

\vspace{0.1in}

%We will prove a super version of the constructive part of this result in
%\S\ref{GeneralizingSchiffmann}.

\section{Dynamical \textit{r}-matrices in the Super Setting}
\label{SectionSuperDynamicalrmatrices}

We can now
concentrate on the super versions of the
above constructions. 

\subsection{Definitions}
\label{DefinitionsSuperDynamical}

Let $\g$ be a simple Lie superalgebra with nondegenerate Killing
form ${(\cdot \;,\cdot)}$. Let ${\h \subset \g}$ be a Cartan
subsuperalgebra, and let $\Delta \subset \h^*$ be the set of
roots associated to $\h$. Fix a set of simple roots $\Gamma$ or
equivalently a Borel $\b$. The \emph{classical dynamical
Yang-Baxter equation} for a meromorphic function ${r:\h^*
\rightarrow \g \otimes \g}$ will be: 
\begin{equation}
\label{SDYBequation}
Alt_s(dr) + 
[r^{12},r^{13}] + [r^{12},r^{23}] + [r^{13},r^{23}] = 0
\end{equation}
\noindent
The differential of $r$ will be defined as
above as:
$$ \begin{matrix}
dr &:& \h^* &\longrightarrow& \g \otimes \g \otimes \g  \\
&& \lambda &\longmapsto& \sum_i x_i \otimes \frac{\partial
r}{\partial x_i} (\lambda)  
\end{matrix}. $$
\noindent
Here $\{x_i\}$ is a basis for $\h$ so all $x_i$ are even.
Recall that ${Alt_s : \g^{\otimes 3} \rightarrow
\g^{\otimes 3}}$ is given on homogeneous elements by: 
$$ Alt_s(a\otimes b \otimes c) = a \otimes b \otimes c  +
(-1)^{|a|(|b|+|c|)} b \otimes c \otimes a + 
(-1)^{|c|(|a|+|b|)} c \otimes a \otimes b, $$
\noindent
In view of all this, we can see that, for $r = \sum_i
{R_i}_{(1)} \otimes {R_i}_{(2)}$:
\begin{align*}
Alt_s(dr) =& 
\sum_i x_i^{(1)} \left(\frac{\partial r }{\partial x_i}
\right)^{(23)} +
\sum_i x_i^{(2)} \left(\frac{\partial r }{\partial x_i}
\right)^{(31)}\\ %\\
+& 
\sum_i (-1)^{|{R_i}_{(1)}| |{R_i}_{(2)}|} x_i^{(3)} \left(
\frac{\partial r}{\partial x_i}\right)^{(12)}.
\end{align*}

For the moment\footnote{See \S\ref{GeneralizingSchiffmann} for a broader
definition.}, we will say that a meromorphic function ${r : \h^*
\rightarrow \g
\otimes \g}$ is a \emph{super dynamical r-matrix with coupling
constant} $\epsilon$ if it is a solution to Equation
\eqref{SDYBequation} and satisfies the generalized unitarity
condition:

\begin{equation}
\label{superdynamicalunitarity} 
r(\lambda) + T_s(r)(\lambda) = \epsilon \Omega,
\end{equation}
\noindent
where $\Omega$ is the Casimir element, i.e. the element of ${\g
\otimes \g}$ corresponding to the Killing form. Here, $T_s$ is the \emph{super
twist} map, defined on the second tensor power of any given super vector space
by: $T_s(a\otimes b) = (-1)^{|a||b|}b \otimes a$. 

The (constant) standard $r$-matrix and its super twist are two easy examples of  super dynamical $r$-matrices. The (constant) non-standard $r$-matrices can also be viewed as simple cases of super dynamical $r$-matrices. Other simple but non-constant examples can be constructed by the theorems that will be presented in the following subsection. We refer the reader to \cite{Kar3} for more examples and discussion.

\subsection{Super Dynamical r-matrices with Zero Weight}
\label{SuperZeroWeight}

A super dynamical r-matrix $r$ is said to satisfy the \emph{zero weight
condition} if:
$$ [h\otimes 1 + 1 \otimes h, r(\lambda)] = 0 \textmd{ for all }
h \in \h, \lambda \in \h^*.$$

In \cite{Kar3}, the author proved the super versions of Theorem 3.2 and Theorem 3.10 from
\cite{EV}, which are, respectively, The Etingof-Varchenko Theorem
\ref{EVtheorem1} and Etingof -Varchenko Theorem \ref{EVtheorem3} of Section
\ref{SchiffmannResult}. The proofs in \cite{Kar3} have clear similarities
to the respective proofs in \cite{EV}. These results basically extend the full classification results of dynamical $r$-matrices with zero weight to the super case. 

Here are the two theorems:

\begin{theorem}
\label{0couple0weighttheorem}
(1) Let $X$ be a subset of the set of roots $\Delta$ of a
simple Lie superalgebra $\g$ with nondegenerate Killing form ${(
\cdot \; , \cdot )}$ such that:

(a) If $\alpha, \beta \in X$ and $\alpha + \beta$ is a root,
then $\alpha + \beta \in X$, and

(b) If $\alpha \in X$, then $-\alpha \in X.$

\noindent
Let $\nu \in \h^*$, and let $D = \sum_{i<j} D_{ij} dx_i \wedge
dx_j$ be a closed meromorphic $2-$form on $\h^*$. If we set
$D_{ij} = -D_{ji}$ for $i \ge j$, then the meromorphic
function:
$$ r(\lambda) =
\sum_{i,j=1}^N D_{ij}(\lambda) x_i \otimes x_j  +  \sum_{\alpha
\in X} \frac{(-1)^{|\alpha|}(e_{\alpha},e_{-\alpha})}{(\alpha,
\lambda-\nu)} e_{\alpha} \otimes e_{-\alpha} $$ 
\noindent
is a super dynamical r-matrix with zero weight and zero
coupling constant.

\noindent
(2) Any super dynamical r-matrix with zero weight and zero
coupling constant is of this form. 
\end{theorem}

\begin{theorem}
\label{0weighttheorem}
(1) Let $\g$ be a simple Lie superalgebra with nondegenerate
Killing form ${(\cdot \; , \cdot )}$. Let $X$ be a subset of the set of
roots $\Delta$ of $\g$ satisfying conditions $(a)$ and $(b)$ of
Theorem \ref{0couple0weighttheorem}.
Pick $\nu \in \h^*$, and define:
$$ \phi_{\alpha} = \left\{ \begin{matrix}
\left( \epsilon / 2\right) \coth\left(
(-1)^{|\alpha|}(e_{\alpha},e_{-\alpha})
\left(\epsilon / 2\right) (\alpha,\lambda - \nu) \right) 
& \textmd{   if  } \alpha \in X \\
\left(\pm\epsilon / 2\right) 
& \textmd{   if  } \alpha \not \in X, 
\textmd{  negative} \\ 
\mp(-1)^{|\alpha|}\left(\epsilon / 2 \right) 
& \textmd{   if  } \alpha \not \in X,
\textmd{  positive}  
\end{matrix} \right. $$
\noindent
Let $D = \sum_{i<j} D_{ij} dx_i \wedge
dx_j$ be a closed meromorphic $2-$form on $\h^*$. If we set
$D_{ij} = -D_{ji}$ for $i \ge j$, then the meromorphic
function:    
$$ r(\lambda) =
\sum_{i,j=1}^N D_{ij}(\lambda) x_i \otimes x_j +
\frac{\epsilon}{2}\Omega + \sum_{\alpha \in \Delta}
\phi_{\alpha}
e_{\alpha} \otimes e_{-\alpha} $$  
\noindent
is a super dynamical r-matrix with zero weight and nonzero
coupling constant $\epsilon$. 

\noindent
(2) Any super dynamical r-matrix with zero weight and nonzero
coupling constant $\epsilon$ is of this form. 
\end{theorem}

\begin{remark}
Note that, if we take the limit as
$\epsilon \rightarrow 0,$ then the statement of Theorem \ref{0weighttheorem} reduces to the
statement of Theorem \ref{0couple0weighttheorem}.
\end{remark}

\subsection{Generalizing Schiffmann's Classification}
\label{GeneralizingSchiffmann}

Here we discuss an extension of the full
classification result of
Schiffmann. The material in this subsection is new, although it is natural and
expected after
\cite{Kar3}. We start with the necessary
terminology, analogous to the non-graded case.

As in the beginning of Section \ref{SectionSuperDynamicalrmatrices}, 
let $\g$ be a simple Lie superalgebra with nondegenerate Killing form
${(\cdot,\cdot)}$. Let $\mathfrak{l} \subset \g$ be a commutative
subsuperalgebra, $\h \subset \g$ the Cartan subsuperalgebra containing
$\mathfrak{l}$, $\Delta \subset \h^*$ the set of roots associated to $\h$,
$\Gamma$ a fixed set of simple roots and $\b$ the associated Borel
subsuperalgebra. 
In this setting we can generalize our definition of super dynamical
$r$-matrices to include all meromorphic solutions of Equation
\eqref{SDYBequation} where the differential $dr$ is given by:
$$ \begin{matrix}
dr &:& \mathfrak{l}^* &\longrightarrow& \g \otimes \g \otimes \g  \\
&& \lambda &\longmapsto& \sum_i x_i \otimes \frac{\partial
r}{\partial x_i} (\lambda)  
\end{matrix}, $$
\noindent
where, this time, $\{x_i\}$ is a basis for $\mathfrak{l} \subset \h$. More
precisely, we will say that a meromorphic function ${r : \mathfrak{l}^*
\rightarrow \g
\otimes \g}$ is a \emph{super dynamical r-matrix with coupling
constant} $\epsilon$ if it is a solution to Equation
\eqref{SDYBequation} and satisfies the generalized unitarity
condition given by Equation \eqref{superdynamicalunitarity}.

A super dynamical r-matrix $r : \mathfrak{l}^* \rightarrow \g \otimes \g$ is
said to be $\mathfrak{l}-$invariant if:
$$ [l\otimes 1 + 1 \otimes l, r(\lambda)] = 0 \textmd{ for all }
l \in \mathfrak{l}, \lambda \in \mathfrak{l}^*.$$
\noindent
The \emph{zero weight condition} of Subsection
\ref{SuperZeroWeight} is easily seen to be equivalent to $\h-$invariance in our
new terminology.

We will say that two super dynamical $r$-matrices are \emph{gauge-equivalent}
if they can be obtained from one another via \emph{gauge-transformations}, see
\S\ref{SchiffmannResult}. 

Let us start with a result about gauge equivalence classes. The non-graded 
version of the following lemma is proved as the first part of the main
classification theorem in Schiffmann's work, \cite{Schi}. The same proof will
work here with no changes, as only the even component of $\g$ is involved 
when determining the gauge transformations to be used.

\begin{lemma}
\label{supergaugetheorem}
Any $\mathfrak{l}-$invariant super dynamical $r$-matrix with coupling constant
$1$ is gauge-equivalent 
to an $\mathfrak{l}-$invariant super dynamical $r$-matrix $\overline{r} :
\mathfrak{l}^* \rightarrow \g \otimes \g$ satisfying:
\begin{equation}
\label{superlskewsymmetry}
\overline{r}(\lambda) - T_s(\overline{r})(\lambda) \in \mathfrak{l}^{\perp}
\otimes 
\mathfrak{l}^{\perp} = \left( (\bigoplus_{\alpha \neq 0}
\g_{\alpha} ) \oplus \h_0 \right) \otimes
\left((\bigoplus_{\alpha \neq 0} \g_{\alpha}) \oplus
\h_0\right)
\end{equation}
\noindent
where $\h_0 \subset \h$ is ``the" complement of $\mathfrak{l}$ in $\h$. $T_s$,
as before, is the super twist, mapping any homogeneous $a \otimes b$ to
$(-1)^{|a||b|}b \otimes a$.
\end{lemma}

\begin{remark}
The orthogonal complement of a subset $\mathfrak{l}$ of the Cartan subalgebra
of a Lie superalgebra may or may not intersect the subset $\mathfrak{l}$ 
trivially. However, after making certain choices, we can always find a 
subset $\h_0 \subset \h$ such that $\h_0 \oplus \mathfrak{l} = \h$.
More specifically, we pick a set $\Gamma$ of simple roots for $\g$ and an
appropriate basis $\{h_{\alpha} \in \h : \alpha \in \Gamma\}$ for $\h$ such 
that for some subset $A$ of $\Gamma$ we
have:
\[ \mathfrak{l} = \bigoplus_{\alpha \in A} \C h_{\alpha}.\]
\noindent
Then we set $\h_0 = \bigoplus_{\alpha \in \Gamma \backslash A} \C h_{\alpha}$.
This is the complement we use in the above statement.
\end{remark}

Let us say that a triple $(\Gamma_1,\Gamma_2,\tau)$ is an \emph{admissible
triple} if:
\begin{enumerate}
\item $\Gamma_1, \Gamma_2 \subset \Gamma,$ and
\item $\tau : \Gamma_1 \rightarrow \Gamma_2$ is a grading preserving isometry.
\end{enumerate}
\noindent
We will say that an admissible triple is $\mathfrak{l}$-\emph{graded} if, in
addition, it preserves
the decomposition of $\g$ into $\mathfrak{l}$-weight spaces. Let $(\Gamma_1,
\Gamma_2, \tau)$ be an
$\mathfrak{l}$-graded admissible triple.

Denote by $\Gamma_3$ the largest subset of the intersection $\Gamma_1 \cap
\Gamma_2$ which is stable under $\tau$, and define: 
\[ \overline{\Gamma}_1 = \Gamma_1 \backslash \Gamma_3, \qquad 
\overline{\Gamma}_2 = \Gamma_2 \backslash \Gamma_3.\]
\noindent
Then, the triple $(\overline{\Gamma}_1, \overline{\Gamma}_2, \tau)$ is
admissible in the sense of \cite{Kar1} and will give us, through the
constructions there, a solution to the CYBE on $\g$.

For each choice of a set of ``Chevalley" generators\footnote{For simple Lie superalgebras with non-degenerate Killing form, it is always possible to find a basis $\{E_{\alpha},F_{\alpha},H_{\alpha}: \alpha \in \Gamma\}$ which satisfies the super versions of the usual commutation and Serre relations one expects from a set of Chevalley generators in the non-graded case. See, for instance, \cite{GL, LS2}.}, we can extend $\tau$ to two (Lie superalgebra) isomorphisms:
\[ \overline{\tau}_{1,2} : \overline{\Gamma}_1 \rightarrow \overline{\Gamma}_2
\quad \textmd{and} \quad
\overline{\tau}_3 : \Gamma_3 \rightarrow \Gamma_3.\]
\noindent
We can and will denote both of these maps by $\overline{\tau}$.

If we fix a basis $\{e_{\alpha} : \alpha \in \Gamma\}$ for the non-Cartan part
of $\g$ and introduce
the $A_{\alpha}$ notation as done before in \S\ref{SuperZeroWeight},
then, for any $\lambda
\in \mathfrak{l}^*$, we can define a map ${K_s(\lambda) :
\mathfrak{n}_+(\Gamma_1)
\rightarrow \mathfrak{n}_+(\Gamma_2)}$ as follows:
\begin{eqnarray*} 
K(\lambda) (e_{\alpha}) &=& \frac{1}{2} e_{\alpha} +
e^{-A_{\alpha}(\alpha,\lambda)}
\frac{\tau}{1-e^{-A_{\alpha}(\alpha,\lambda)}\tau}(e_{\alpha}) \\
&=& \frac{1}{2}e_{\alpha} + 
\sum_{n>0} e^{-nA_{\alpha}(\alpha,\lambda)}\tau^n(e_{\alpha}). 
\end{eqnarray*}
\noindent
Note that the sum above is finite as long as ${\alpha \not\in
\Gamma_3}$. If ${\alpha \in \Gamma_3}$ and $\tau|_{\Gamma_3} = id_{\Gamma_3}$,
then we have:
\[ K(\lambda) (e_{\alpha}) = \frac{1}{2} 
\coth\left(\frac{A_{\alpha}}{2}(\alpha,\lambda) \right)e_{\alpha}.
\]

Let ${r_{00} \in \h_0 \otimes \h_0}$ satisfy: 
\begin{equation*}
(\tau(\alpha) \otimes 1) r_{00}  + (1 \otimes \alpha)r_{00} =
\frac{1}{2} ((\alpha + \tau(\alpha)) \otimes 1) \Omega_{00}
\end{equation*}
\noindent
for any $\alpha \in \Gamma_1$, where we denote by $\Omega_{00}$ the $\h_0
\otimes \h_0$ part of the
Casimir element $\Omega$.
Then some computation shows in fact that 
\[ r(\lambda) = \frac{1}{2}\Omega + r_{00} + \sum_{\alpha \in
\overline{\Gamma}_1 \cap \Delta^+} K(\lambda)(e_{\alpha}) \wedge
e_{-\alpha} + \sum_{\alpha \in \Delta^+, \alpha \not \in
\overline{\Gamma}_1} \frac{1}{2}e_{\alpha} \wedge e_{-\alpha}\]
\noindent
is a dynamical r-matrix satisfying Equation
\eqref{superlskewsymmetry}. 
Summarizing the above, we get:

\begin{theorem}
\label{superschiffthm}
Let ${(\Gamma_1, \Gamma_2, \tau)}$ be an $\mathfrak{l}-$graded admissible
triple and let ${r_{00} \in \h_0 \otimes \h_0}$ satisfy: 
\begin{equation*}
(\tau(\alpha) \otimes 1) r_{00}  + (1 \otimes \alpha)r_{00} =
\tfrac{1}{2} ((\alpha + \tau(\alpha)) \otimes 1) \Omega_{00}
\end{equation*}
\noindent
for any $\alpha \in \Gamma_1$. Then:
$$ r(\lambda) = \frac{1}{2}\Omega + r_{00} + \sum_{\alpha \in
\overline{\Gamma}_1 \cap \Delta^+} K(\lambda)(e_{\alpha}) \wedge
e_{-\alpha} + \sum_{\alpha \in \Delta^+, \alpha \not \in
\overline{\Gamma}_1} \frac{1}{2}e_{\alpha} \wedge e_{-\alpha}$$
\noindent
is a dynamical r-matrix satisfying Equation \eqref{superlskewsymmetry}.
\end{theorem}

Together with Lemma \ref{supergaugetheorem}, Theorem \ref{superschiffthm} gives
us the super version
of the constructive part of Schiffmann's result. Hence, we have a very nice way
to construct super
dynamical $r$-matrices which generalizes naturally the non-graded theory. Note
that, as in the non-graded case, this result agrees with and extends the constructive part of the zero weight results of \S\ref{SuperZeroWeight}. More precisely, if we let $\mathfrak{l} = \h$, $\Gamma_1 = \Gamma_2 = X$, $\tau = id$ and $\epsilon = 1$, then Theorems \ref{0weighttheorem} and \ref{superschiffthm} coincide.

However, superizing the classifying part of the non-graded theory will be a lot
more involved. The
$sl(2,1)$ example constructed in \cite{Kar1} and studied in detail in
\cite{Kar2} once again turns
out to be an issue. It can easily be seen that for $\Gamma_3 = 0$ and
$\mathfrak{l} = 0$, Theorem
\ref{superschiffthm} reduces to the construction theorem of \cite{Kar1},
and the construction described above gives us a non-skewsymmetric solution to
the CYBE. The
particular $sl(2,1)$ example which does not fit the framework
of \cite{Kar1}
therefore will not fit this new framework.

\section{General Theory of Dynamical Quantum Groups}
\label{SectionDynamicalQuantumGroups}

Here we review the general theory of dynamical quantum groups. In order to keep the paper at a readable length, we only consider the development of the theory utilizing the notion of Hopf algebroids. Nevertheless, we should note that
there are alternative approaches. For instance, one can use quasi-Hopf structures (i.e. Hopf-like structures obtained by weakening the coassociativity condition on
the coproduct), which were first introduced by Drinfeld in \cite{D90}. Such an approach was initiated in \cite{BBB}, and its implications were investigated in great depth; see, for instance, \cite{ABRR, EF}. It turns out that there is a natural relationship between quasi-Hopf algebras and Hopf algebroids (\cite{X99}). In fact, the reader interested in details on using different generalizations of Hopf algebras in the context of dynamical quantum groups may find our survey \cite{Kar4} of some value. However, in this note, we will not be too worried about only using the Hopf algebroid approach without any further explanation.

\subsection{Groupoids, Bialgebroids and Hopf Algebroids}
\label{WhatisHopfAlgebroid}

It is well-known that quantum groups are actually Hopf algebras. It turns out that,
%the correct Hopf objects to consider 
in the context of dynamical quantum groups, 
considering Hopf algebroids as the analogous objects proves quite fruitful.
%a good and fruitful choice of Hopf-like objects is going to be Hopf algebroids. 
Therefore, we will start with a
basic discussion of groupoids, bialgebroids, Hopf algebroids, and quantum
groupoids.

For our purposes, a \emph{groupoid over a set $X$} is a set $G$ together with the following structure maps:
\begin{enumerate}
\item A pair of maps $s, t : G \rightarrow X$, respectively called the
\emph{source} and the 
\emph{target}.
\item A \emph{product} $m$, i.e. a partial function $m : G \times G \rightarrow
G$ satisfying the
following two properties:
\begin{enumerate}
\item $t(m(g,h)) = t(g)$, $s(m(g,h)) = s(h)$ whenever $m(g,h)$ is defined;
\item $m$ is associative: $m(m(g,h),k) = m(g,m(h,k))$ whenever the relevant
terms are defined. 
\end{enumerate}
\item An embedding $\epsilon : X \rightarrow G$ called the \emph{identity
section} such that
$m(\epsilon(t(g)),g) = g = m(g, \epsilon(s(g)))$ for all $g \in G$.
\item An \emph{inversion map} $\imath : G \rightarrow G$ such that
$m(\imath(g),g) =
\epsilon(s(g))$ and $m(g,\imath(g)) = \epsilon(t(g))$ for all $g \in G$.
\end{enumerate}

Conventionally one writes $gh$ instead of $m(g,h)$ whenever the latter is
defined. Also, $\imath(g)$
can be denoted by $g^{-1}$. Then we can rewrite the above conditions in a form
which makes
more transparent the similarities and differences of \emph{a groupoid}
from \emph{a group}. For
instance, we wish to have inverses for all elements of $G$, but the
multiplication is only partially
defined.

There are other ways to define a groupoid. One very elegant way is to view it as a
particular type of category $(X,G)$ with $X$ making up the set of objects, such that the morphisms (elements of $G$) are all invertible. In this more category-theoretic setup, the notion of a group corresponds to the particular case when X has only one object, and this single object is the domain and range for all morphisms (the group elements), just like the set of morphisms of a category with only one object corresponds, in general, to the concept of a monoid. We choose not to consider this categorical description. However, if we draw some schematic figures to represent the structure maps defined above, then we can clearly see how the categorical notion may be derived easily. This may also help interpret, in graph-theoretical terms, the particular names \emph{source} and \emph{target} used for the two structure maps in the above definition.

Groupoids were first introduced in 1926, and since then, found
applications in differential topology and geometry, algebraic geometry and
algebraic topology, and analysis. A very friendly introduction to groupoids
with many examples from various areas of mathematics can be found in
\cite{Wei2}. For more rigorous accounts one may refer to the bibliography
there. We took our definition from \cite[Part \textrm{VI}]{CW}. Another good reference on groupoids is \cite{Pat99}. 

It is possible to define Lie groupoids or Poisson groupoids by stipulating the
presence of a Lie or Poisson structure on a given groupoid, but since we will not need precise definitions here, we will only mention a few references for the interested reader. For instance, \cite{CDW} provides a systematic introduction to Lie groupoids, Lie algebroids and symplectic groupoids. Poisson groupoids were first introduced in \cite{Wei}, and they have been studied ever since. The recent monograph \cite{LP} on dynamical Poisson groupoids studies in great detail the particular Poisson groupoids most relevant to this note.
 
A possible connection of Poisson groupoids to quantum groupoids analogous to
the connection of Poisson groups to quantum groups was conjectured in 
\cite{May}.\footnote{This reference also includes a brief but interesting philosophical discussion of quantization.} The first full description of a
quantum groupoid as an object which should simultaneously
be generalizing quantum groups and groupoid algebras was given in
\cite{Mal}. This latter work utilized certain commutativity assumptions, and
eventually a new and broader approach, in \cite{Lu1}, provided us with a
description without such constraints.

Almost simultaneously, the notion of a Hopf algebroid was being developed.
Commutative Hopf algebroids were first studied in \cite{Rav}. There were
similar descriptions in \cite{Mal}. In \cite{Lu1}, the definitions were
slightly modified in order to include noncommutative cases. However, the clear
consensus on the definition of a bialgebroid, first given in \cite{Tak}, was
not easy to come for the definition of a Hopf algebroid. The antipode
suggested in \cite{Lu1} was not universally accepted, and various other
formulations followed. See, for instance, \cite{BS2}. A comparative study of
these various antipodes may be found in \cite{Bo}. 

Another parallel development was the introduction of weak Hopf algebras in
\cite{BNS, BS1, N}. These are the Hopf-algebra-like structures one obtains when
one drops the requirement that the comultiplication be unit preserving. They
were introduced with a view toward applications to operator algebras, but even
from the beginning, their appeal as a means of generalizing quantum groups was
recognized. Eventually, B\"{o}hm and Szlach\'{a}nyi showed that weak Hopf
algebras with bijective antipodes are Hopf algebroids, see \cite{BS2}. The
reader can also refer to \cite{NV} for a detailed overview of weak Hopf
algebras, their relationship to various generalizations of the idea of quantum
groups and in particular to ``\emph{dynamical deformations of quantum
groups}." The survey \cite{Kar4} provides a basic comparative study, in the context of possible generalizations of Hopf algebras. 

Here are the definitions we use in this paper, mainly following
\cite{Bo} and \cite{BS2}:
A bialgebroid should be the natural extension of the notion of a bialgebra to
the world of groupoids. This will imply that a bialgebroid is no longer an
algebra, but a bimodule over a non-commutative ring. More specifically, a
\emph{left bialgebroid} $\mathcal{A}_L$ is given by the following data: 
\begin{enumerate}
\item Two associative unital rings: the \emph{total ring} $A$ and the
\emph{base field} $L$.
\item Two ring homomorphisms: the \emph{source} ${s_L : L \rightarrow A}$ and
the \emph{target} ${t_L : L^{op} \rightarrow A}$ such that the images of $L$
in $A$ commute, making $A$ an $L-L$ bimodule denoted by $_LA_L$.
\item Two maps ${\gamma_L : A \rightarrow A_L \otimes {_LA}}$ and ${\pi_L : A
\rightarrow L}$ making the triple ${(_LA_L, \gamma_L, \pi_L)}$ a comonoid in
the category of $L-L$ bimodules. 
\end{enumerate}

The source and target maps $s_L$ and $t_L$ may be used to define four commuting
actions of $L$ on $A$; these in turn give us in an obvious way the new
bimodules $^LA^L$, $^LA_L$, and $_LA^L$.

Similarly we can define a \emph{right bialgebroid} $\mathcal{A}_R$ using the
following data:
\begin{enumerate}
\item Two associative unital rings: the \emph{total ring} $A$ and the
\emph{base field} $R$.
\item Two ring homomorphisms: the \emph{source} ${s_R : R \rightarrow A}$ and
the \emph{target} ${t_R : R^{op} \rightarrow A}$ such that the images of $R$
in $A$ commute, making $A$ an $R-R$ bimodule denoted by $^RA^R$.
\item Two maps ${\gamma_R : A \rightarrow A^R \otimes {^RA}}$ and ${\pi_R : A
\rightarrow R}$ making the triple ${(^RA^R, \gamma_R, \pi_R)}$ a comonoid in
the category of $R-R$ bimodules. 
\end{enumerate}

As in the case of left bialgebroids, we can define three other bimodule
structures on $A$ using the source and the target, and denote them by $_RA_R$,
$^RA_R$, and $_RA^R$. These bimodule structures and the two notions of
bialgebroids are related as expected. For instance, if $\mathcal{A}_L = (A, L,
s_L, t_L, \gamma_L, \pi_L)$ is a left bialgebroid, then its co-opposite is
again a left bialgebroid: $(\mathcal{A}_L)_{cop} = (A, L^{op}, t_L, s_L,
\gamma_L^{op}, \pi_L)$, where $\gamma_L^{op}$ is defined as $T \circ
\gamma_L$.\footnote{Here, as before, $T$ is the usual (non-graded) twist, mapping $a \otimes b$ to $b \otimes a$.} The opposite $(\mathcal{A}_L)^{op}$ defined by the data  $(A^{op}, L, t_L, s_L, \gamma_L, \pi_L)$ is a right bialgebroid. For more on bialgebroids, we refer the reader to \cite{Bo}.

We will take our definition for a Hopf algebroid from \cite{BS2}. In
particular, to define a Hopf algebroid, we will need two associative unital
rings $A$ and $L$, and set $R = L^{op}$. We will consider a left bialgebroid
structure $\mathcal{A}_L = (A, L, s_L, t_L, \gamma_L, \pi_L)$ and a right
bialgebroid structure $\mathcal{A}_R = (A, R, s_R, t_R, \gamma_R, \pi_R)$
associated to this pair of rings. We will require that $s_L(L) = t_R(R)$ and
$t_L(L) = s_R(R)$ as subrings of $A$, and:
\begin{eqnarray*}
(\gamma_L \otimes id_A) \circ \gamma_R &=& (id_A \otimes \gamma_R) \circ
\gamma_L, \\
(\gamma_R \otimes id_A) \circ \gamma_L &=& (id_A \otimes \gamma_L) \circ
\gamma_R. 
\end{eqnarray*}
The last ingredient is the antipode. This will be a bijection $S : A
\rightarrow A$ that will satisfy:
\begin{eqnarray*}
S(t_L(l)at_L(l')) &=& s_L(l')S(a)s_L(l), \\
S(t_R(r')at_R(r)) &=& s_R(r)S(a)s_R(r')
\end{eqnarray*}
\noindent
for all $l,l' \in L$, $r,r' \in R$, and $a \in A$. In other words, we require
$S$ to be a twisted isomorphism simultaneously of bimodules $^LA_L \rightarrow
{_LA^L}$ and of bimodules $^RA_R \rightarrow {_RA^R}$. 

Our final constraint on the antipode $S$ is as follows: 
\begin{eqnarray*}
S(a_{(1)})a_{(2)} &=& s_R \circ \pi_R(a), \\
a^{(1)}S(a^{(2)}) &=& s_L \circ \pi_L(a)
\end{eqnarray*}
\noindent
for any $a \in A$. The subscripts and the superscripts on $a$ come from a
generalized version of the famous Sweedler notation, which we use to define the
two maps $\gamma_L$ and $\gamma_R$: 
\begin{eqnarray*}
\gamma_L(a) &= a_{(1)} \otimes a_{(2)} &\in A_L \otimes {_LA} \\
\gamma_R(a) &= a^{(1)} \otimes a^{(2)} &\in A^R \otimes {^RA}.
\end{eqnarray*}

In this setup, then, we will say that the triple $ \mathcal{A} =
(\mathcal{A}_L, \mathcal{A}_R, S)$ is a \emph{Hopf algebroid}.\footnote{We should note here that some of the above information in our definition is redundant. In fact one can start with a left bialgebroid $\mathcal{A}_L = (A, L, s_L, t_L, \gamma_L, \pi_L)$ and an anti-isomorphism $S$ of the total ring $A$ satisfying certain conditions, and from here can reconstruct a right bialgebroid $\mathcal{A}_R$ using the same total ring such that the triple $(\mathcal{A}_L, \mathcal{A}_R, S)$ is a Hopf algebroid. See \cite{Bo} for more details.}
It is easy to see how this symmetric definition, in terms of two bialgebroids
and a bijection called the \emph{antipode}, is analogous to the definition of
a Hopf algebra from two bialgebras and a bijective map called an antipode.
However, it is not nearly as easy to see why this is the appropriate
definition. We accept the definition given above without further
analysis, and leave the readers to follow the discussion on this issue on
their own. (A good place to start may be the comparative study of B\"{o}hm in \cite{Bo}).

Now, for us, a quantum groupoid will be a particular type of Hopf algebroid,
just as a quantum group is a particular type of Hopf algebra. Quantum groups
are Hopf algebras obtained from deformations of commutative or cocommutative
Hopf algebras. The semiclassical limits of these deformations are the
so-called Hopf-Poisson algebras; associated to these latter structures are the
Poisson-Lie groups. Hence, a quantum groupoid should be a Hopf algebroid which
is a ``deformation" of a ``nice Hopf algebroid" in such a way that in the
semi-classical limit we should get an algebraic structure associated to a
Poisson-Lie groupoid.

We will end this section with the rather intuitive discussion above, 
and not attempt to come up with a full accurate definition for the notion of
\emph{quantum groupoid}. Our main reason for this is the fact that just as
there are various definitions for Hopf algebroids, there are some differing
notions of quantum groupoids. For instance, in \cite{NV}, the term \emph{quantum groupoid} is
used almost interchangeably with the term \emph{weak Hopf algebra}. (Also see
\cite{X01} for another different approach to quantum groupoids).
However, unlike in the case of
Hopf algebroids, we can do without a precise definition at this stage, because
there is nevertheless a consensus on the definition of dynamical quantum
groups\footnote{These structures should technically be called \emph{dynamical
quantum groupoids}, as they are certain types of quantum groupoids related to
the quantum dynamical Yang-Baxter equation (QDYBE). However, the term
\emph{dynamical quantum group} is more common in the literature.}. We will describe and study those in more detail in the next subsection.

\subsection{Dynamical Quantum Groups}
\label{QuantumGroupoids} 

Here we summarize the current theory of dynamical quantum groups.
We begin with the basic definitions and then provide the necessary
connections with the algebraic terms from the previous subsection. We mainly
follow the theory as developed in \cite{EV2} and summarized in \cite{ES2}, and
explain how the earlier discussions of \S\ref{WhatisHopfAlgebroid}
are compatible with it.\footnote{At this point, we should remark that the antipode as defined in \cite{EV2} contains small inconsistencies, which were noted and modified in \cite{KR}. Therefore, it is more accurate to say that we will follow \cite{EV2} up to some corrections made in \cite{KR}. This comment will become clearer in the following discussion.}

In order to define dynamical quantum groups, we start with the notion of an
\emph{$H$-algebra}, where $H$ is taken to be a commutative, cocommutative,
finitely generated Hopf algebra over $\C$. Let $G = \Spec H$ be the
corresponding commutative affine algebraic group. Assuming that $G$ is
connected, let $M_G$ denote the field of meromorphic functions on $G$. Then we
will say that an associative unital $\C$-algebra $A$ is an \emph{$H$-algebra}
if it has a $G$-bigrading $A = \oplus_{\alpha,\beta \in G} A_{\alpha,\beta}$
(called the \emph{weight decomposition}), and two algebra embeddings $\mu_l ,
\mu_r : M_G \rightarrow A_{0,0}$ (called the \emph{left} and \emph{right
moment maps}\footnote{These structure maps generalize certain maps which are interchangeably called \textit{moment maps}, \textit{momentum maps}, or \textit{momentum mappings} in the literature. For purely typographical reasons, we will be using the shortest phrase among these three candidates.} respectively) such that we have:
\begin{eqnarray*}
\mu_l(f(\lambda))a &=& a \mu_l(f(\lambda + \alpha)), \\
\mu_r(f(\lambda))a &=& a \mu_r(f(\lambda + \beta)),
\end{eqnarray*}
\noindent
for any $a \in A_{\alpha,\beta}$ and $f \in M_G$. With a slight change of
perspective, we can summarize the above by saying that the ordered quadruple
$(A, H, \mu_l, \mu_r)$ is an $H$-algebra.

Next we define the algebra $D_G$ of difference operators on $M_G$. In other
words, $D_G$ consists of all operators ${M_G \rightarrow M_G}$ of the form
$\sum_{i=1}^n f_i(\lambda)\sigma_{\beta_i}$ where $f_i \in M_G$ and for
${\beta \in G}$, $\sigma_{\beta}$ is the automorphism of $M_G$ given by
$\sigma_{\beta}(f) (\lambda) = f(\lambda + \beta)$. 
Note that, using these $\sigma_{\alpha}$, $\sigma_{\beta}$, we can rewrite the
conditions on the left and right moment maps given above as follows:
\[  \mu_l(f)a = a \mu_l(\sigma_{\alpha}f), \qquad
\mu_r(f)a = a \mu_r(\sigma_{\beta}f) \]
\noindent
Now, $D_G$ itself is an $H$-algebra where the bigrading is given by:
$f\sigma_{-\alpha} \in (D_G)_{\alpha,\alpha}$. Both moment maps are taken to
be the natural embedding of $M_G$ in $D_G$. 
In fact, it can be shown that $D_G$ behaves very much like the unit object in
the category of $H$-algebras in the following sense: Given any $H$-algebra
$A$, there are canonical $H$-algebra isomorphisms $A \cong A \otimes D_G 
\cong D_G \otimes A$ where $a \in A_{\alpha,\beta}$ is mapped to $a \otimes
\sigma_{-\beta}$ and $\sigma_{-\alpha} \otimes a$, respectively, see \cite{KR}
for details. 

We can now give the definition for an \emph{$H$-bialgebroid}. 

For this, we need a coassociative \emph{coproduct}, that is, a homomorphism of $H$-algebras ${\Delta : A \rightarrow A \otimes A}$ which satisfies:
\[{(\Delta \otimes id_A) \circ \Delta = (id_A \otimes
\Delta) \circ \Delta.}\] 
\noindent
We also need a \emph{counit}, that is, a
homomorphism of $H$-algebras ${\epsilon : A \rightarrow D_G}$ which satisfies the \emph{counit axiom}: 
\[ (\epsilon \otimes id_A) \circ \Delta = (id_A \otimes \epsilon) \circ \Delta
= id_A. \]
\noindent
In short, we say that the ordered quintuple $\mathcal{A}_H = (A, H, \mu_l, \mu_r ,\Delta,
\epsilon)$ is an \emph{$H$-bialgebroid} if $(A, H, \mu_l, \mu_r)$ is an
$H$-algebra, and $\Delta$ and $\epsilon$ satisfy the required conditions
above. 

To obtain an $H$-Hopf algebroid, we only need to add to the bag an
\emph{$H$-antipode}, i.e. an antiautomorphism ${S : A \rightarrow A}$
satisfying the following technical condition:
\begin{eqnarray*}
a_{(1)}S(a_{(2)}) &=& \mu_l(\epsilon(a)1), \\
S(a_{(1)})a_{(2)} &=& \mu_r(\sigma_{\alpha}(\epsilon(a)1)),
\end{eqnarray*}
\noindent
where $a \in A_{\alpha,\beta}$ and we used the Sweedler notation
associated to $\Delta(a)$. $\epsilon(a)1$ is the result of applying the
difference operator $\epsilon(a) \in D_G$ to the constant function $1 \in
M_G$.\footnote{Note that since $\epsilon$ is an $H$-algebra homomorphism, it
preserves the bigrading and so $\epsilon(A_{\alpha,\beta}) = 0$ unless $\alpha
= \beta$. Therefore, in the second equation, we could have written
$\sigma_{\beta}$ instead of $\sigma_{\alpha}$.}
With this last ingredient we can finally define an
\emph{$H$-Hopf algebroid}: It is an $H$-bialgebroid with an
$H$-antipode.\footnote{In \cite{EV2}, when the antipode for an $H$-bialgebroid is
defined, the term $\sigma_{\alpha}$ is missing from the second equation.
However, the correction term is necessary for making the rest of the arguments
follow consistently; see \cite{KR} for more details. Also, it may be
interesting to note that the authors of \cite{KR} choose to make a much weaker
definition, and then, they prove that one can still get all of the desired properties of an antipode (e.g. uniqueness, antiautomorphism).}

At this point, it is natural to ask how these structures relate to the
structures we defined in the previous subsection. A brief comparison will in
fact show that an $H$-bialgebroid as defined in this subsection is going to be
a (left) bialgebroid in the sense of \S\ref{WhatisHopfAlgebroid}.
More specifically, an $H$-bialgebroid $\mathcal{A}_H = (A, H, \mu_l, \mu_r,
\Delta, \epsilon)$ is a special type of a (left) bialgebroid where the total
ring is $A$ and the base field is $L = M_G$. The source and the target maps
for the left bialgebroid are determined by the two quantum moment maps $\mu_l$
and $\mu_r$. 
The coproduct $\Delta$ and the counit $\epsilon$ are, respectively, the two
maps $\gamma_L$ and $\pi_L$ making $A$ into a comonoid in the category of
$L-L$ bimodules. 
Similarly adding an $H$-antipode to an $H$-bialgebroid to obtain
an $H$-Hopf algebroid is seen to be equivalent to the addition of the antipode
to get a Hopf algebroid in \S\ref{WhatisHopfAlgebroid}. The $\sigma_{\alpha}$
will pop up in the third antipode equation while we use the $H$-structure to
define a right bialgebroid structure.

How do $H$-Hopf algebroids come up in the realm of quantum groups? To answer this question, we first need to talk about the quantum dynamical Yang-Baxter equation (QDYBE). 

Let $\h$ be a finite-dimensional abelian Lie algebra, and let $V$ be a semisimple $\h$-module. Then, for a meromorphic function $R
: \h^* \rightarrow End_{\h}(V \otimes V)$, the 
%\emph{quantum dynamical Yang-Baxter equation} 
QDYBE is the following equation in $V \otimes V \otimes V$:
\begin{equation}
\label{QDYBE}
R^{12}(\lambda- h^{(3)}) R^{13}(\lambda) R^{23}(\lambda-h^{(1)}) =
R^{23}(\lambda) R^{13}(\lambda - h^{(2)}) R^{12}(\lambda).
\end{equation}
\noindent 
Here, $h^{(i)}$, $i = 1,2,3$, is to be replaced by $\mu$ if $\mu$ is the weight
of the $i$th tensor component. An invertible solution $R$ of the QDYBE is
called a \emph{(quantum) dynamical $R$-matrix}. Now, certain ``nice" dynamical
$R$-matrices, in particular the ones satisfying the so-called \emph{Hecke
condition}, may be used to define $H$-Hopf algebroids, where $H$ is the universal enveloping algebra of $\h$. 

A dynamical $R$-matrix $R$ is said to satisfy the \emph{Hecke condition}, with parameter $q \in \C^*$ if the eigenvalues of $TR$ are $1$ on the weight subspaces of $V \otimes V$ of the type $V_a \otimes V_a$ and $1, q$ on the weight subspaces of the type $(V_a \otimes
V_b) \oplus (V_b \otimes V_a)$. As noted in \cite{ES2}, the Hecke condition is a quantum version of Equation \eqref{dynamicalunitarity}, the generalized unitarity condition. In
particular, if a continuous family of dynamical $R$-matrices of the form
$R = 1 - \gamma r + O(\gamma^2)$ satisfies the Hecke condition with
parameter $q =1$, then the semiclassical limit $r$ satisfies the
unitarity condition $r + T(r) = 0$.

We do not go into more
detail here and refer the reader to \cite{EV2, ES2} for more information. Also
see the next section where we discuss superized versions of the QDYBE,
dynamical $R$-matrices and the Hecke condition.

Finally we are in a position to define \emph{dynamical quantum groups}: These
are going to be those $H$-Hopf algebroids which can be obtained from dynamical
$R$-matrices satisfying the Hecke condition. Clearly they are Hopf algebroids
in the sense of \S\ref{WhatisHopfAlgebroid} and are related to the
QDYBE, as expected. 

\subsection{The Categorical Picture}
\label{CategoricalSubsection}

To understand some of the constructions above, and to be able to extend them to the super case, it is imperative that we spend some time on their category-theoretic foundations. The relevant framework was first developed in \cite{EV2}. In our presentation, we will follow \cite{ES2}.

Recall that an algebra $A$ that is also a coalgebra is a bialgebra precisely when the category of modules over $A$ is a monoidal category with the action of $A$ on any tensor product of modules being induced by the comultiplication. In this sense, bialgebras are algebraic counterparts of monoidal categories. In a similar sense, Hopf algebras are the algebraic counterparts of rigid monoidal categories.
In order to talk about quasitriangular structures, $R$-matrices and quantum
groups, we need to think further of braided monoidal categories. 

More
generally, if $\mathcal{B}$ is a braided monoidal category, $V$ a symmetric
tensor category 
%(e.g. the category of complex vector spaces, representations
%of a bialgebra or a Hopf algebra) 
and ${\mathcal{F} : \mathcal{B} \rightarrow
V}$ a tensor functor, then using any object $X$ of $\mathcal{B}$, we can
construct a $V$-automorphism of $\mathcal{F}(X) \otimes \mathcal{F}(X)$ which
satisfies the QYBE.\footnote{Some readers may be more familiar with the braid equation, which is a close relative of the QYBE. It was indeed the connection between these two equations that allowed the construction of link invariants from the theory of quantum groups.}
%quantum Yang-Baxter equation. 
See, for example, \cite{CP, ES} for more on monoidal and braided monoidal categories and their relevance to the theory of quantum groups.

In order to talk about the dynamical Yang-Baxter equation and its solutions, we
will again need a braided monoidal category $\mathcal{B}$ and a tensor category
$V$. However, this time $V$ will not typically be a category of modules over a
bialgebra or Hopf algebra. Instead, we will consider the so-called
\emph{dynamical representations} of $H$-bialgebroids and $H$-Hopf algebroids.
Below is the definition for this special type of representation, following
\cite{ES2}:

A \emph{dynamical representation} of an $H$-algebra $A$ is a pair $(W,\pi_W)$,
where $W$ is a diagonalizable $H$-module and $\pi_W : A \rightarrow D_{G,W}$ is
a homomorphism of $H$-algebras. Here, as in \S\ref{QuantumGroupoids},
$H$ is taken to be a commutative, cocommutative, finitely generated Hopf
algebra over $\C$, $G = \Spec H$ is the corresponding commutative affine
algebraic group (assumed to be connected), $M_G$ is the field of meromorphic
functions on $G$, and $D_G$ is the algebra of difference operators on $M_G$.
In this setup, we define $D_{G,W}$ to be the algebra of all difference
operators on $G$ with coefficients in $End_{\C}(W)$: ${D_{G,W} =
\oplus_{\alpha} D^{\alpha}_{G,W}}$, where ${D^{\alpha}_{G,W} \subset
Hom_{\C}(W,W\otimes D_G)}$ is the space of all difference operators on $G$ with
coefficients in $End_{\C}(W)$ that have weight $\alpha \in G$ with respect to
the action of $H$ in $W$.

Next, we need to define tensor products of dynamical
representations in order to obtain the right kind of category. For instance, following
\cite{ES2, EV2}, we let the
tensor product $V \hat{\otimes} W$ of two dynamical representations $V$ and
$W$ of an $H$-bialgebroid $A$ to be the usual tensor product $V \otimes W$ in
the category of vector spaces, and define the map $\pi_{V
\hat{\otimes} W} : A \rightarrow D_{G,V \hat{\otimes} W}$ by setting $\pi_{V
\hat{\otimes} W}(a) = \theta_{V,W} \circ (\pi_V \otimes \pi_W) \circ
\Delta(a)$. Here $\theta_{V,W}$ is the natural embedding, as an
$H$-algebra, of $D_{G,V} \hat{\otimes} D_{G,W}$ into $D_{G, V
\hat{\otimes} W}$ given by:
\[ f\sigma_{\beta} \otimes g \sigma_{\delta} \longmapsto (f \hat{\otimes} g)
\sigma_{\delta} \]
\noindent
where $\sigma_{\beta}$ is the automorphism of $M_G$ (embedded naturally inside
$D_G$), given by $\sigma_{\beta}(f) (\lambda) = f(\lambda + \beta)$ and
$\sigma_{\delta}$ is defined in a similar manner. The function $f
\hat{\otimes} g$ is determined by: 
\[f \hat{\otimes} g (\lambda)(v \otimes w) = f(\lambda-\mu)(v) \otimes
g(\lambda)(w) \]
\noindent
if $g(\lambda)(w)$ has weight $\mu$. This definition of tensor products makes the category of dynamical representations of an $H$-bialgebroid a tensor category. Following \cite{EV2},  we can also define left and right duals for dynamical representations of $H$-Hopf algebroids, but we will not do so here.

Now if $V_d$ is the tensor category of dynamical representations of an $H$-bialgebroid or an $H$-Hopf algebroid, and if we are given a tensor functor ${\mathcal{F} : \mathcal{B} \rightarrow V_d}$, then using any object $X$ of $\mathcal{B}$, we can once again construct a
$V_d$-automorphism $R = {R(\mathcal{B},\mathcal{F},X)}$ of $\mathcal{F}(X)
\otimes \mathcal{F}(X)$. It turns out that this new automorphism does not
satisfy the quantum Yang-Baxter equation but instead the quantum dynamical
Yang-Baxter equation (QDYBE). Conversely, given a dynamical $R$-matrix $R$,
we can find appropriate categories $\mathcal{B}$ and $V_d$,
a tensor functor ${\mathcal{F} : \mathcal{B} \rightarrow V_d}$ and a
particular object $X$ of $\mathcal{B}$ such that $R$ is the $V_d$-automorphism
associated to the triple ${(\mathcal{B},\mathcal{F},X)}$ in the same manner as
above, see \cite{ES2}.

\section{Dynamical Quantum Groups in the Super Setting}
\label{SectionSuperGroups}

We now begin our study of the dynamical quantum super groups. Once again we
emphasize that we will exclusively follow the Hopf algebroid approach for
dynamical quantum groups. We should note that a superization of the quasi-Hopf
algebra approach for dynamical quantum groups was already initiated, in \cite{GZh, GIZh}.\footnote{Incidentally, an earlier result related to dynamical quantum
super groups which
is definitely worth mentioning can be found in \cite{Is}. Here, the author
works in the framework of the linear quantum groups and computes several
solutions to the QDYBE in both the non-graded and the graded cases. As we
develop our theory based on Hopf structures, it will be interesting to see how
our methods relate to those used in \cite{Is} and compare results. }

\subsection{Dynamical Quantum Groups - The Super Story}
\label{DefinitionsSuperGroups}

In order to get the correct definitions for the super analogues of the notions of Hopf algebroids and dynamical quantum groups, we now concentrate on the structures described in \S\ref{QuantumGroupoids} and \S\ref{CategoricalSubsection}, and superize them systematically. Most of the superization will be straight forward, but we still wish to state our definitions explicitly as often as possible. Some repetition, therefore, will be unavoidable, but hopefully, this will help the reader follow the paper with more ease.

We start with the superization of the constructions in \S\ref{QuantumGroupoids}. In particular, we start with the notion of an \emph{$H$-superalgebra}, where $H$ is a commutative, cocommutative, finitely generated super Hopf algebra over $\C$.\footnote{Super Hopf algebras are special examples of Hopf algebras in the braided monoidal category of Yetter-Drinfeld modules, see \cite{CG, Maj, RT}
for details and important results on Yetter-Drinfeld modules.}
Set $G = \Spec H$, and assuming that $G$ is connected, let $M_G$ denote the field of
meromorphic functions on $G$. Then we will say that an associative unital
$\C$-superalgebra $A$ is an \emph{$H$-superalgebra} if it has a $G$-bigrading
$A = \oplus_{\alpha,\beta \in G} A_{\alpha,\beta}$ (called the \emph{weight
decomposition}) compatible with the $\Z/{2\Z}$-grading\footnote{In other words, we have: $A_{\alpha,\beta} = (A_{\alpha,\beta})_{\overline{0}} \oplus (A_{\alpha,\beta})_{\overline{1}}$, where $(A_{\alpha,\beta})_{\overline{p}} = A_{\alpha,\beta} \cap A_{\overline{p}}$.}, and two superalgebra
embeddings ${\mu_l, \mu_r : M_G \rightarrow A_{0,0}}$ with:

\begin{eqnarray*}
\mu_l(f(\lambda))a &=& a \mu_l(f(\lambda + \alpha)), \\
\mu_r(f(\lambda))a &=& a \mu_r(f(\lambda + \beta)),
\end{eqnarray*}
\noindent
for any $a \in A_{\alpha,\beta}$ and $f \in M_G$.  ($\mu_l, \mu_r$ are called
the \emph{left}
and \emph{right moment maps}). With a slight change of
perspective, we can summarize the above by saying that the ordered quadruple
$(A, H, \mu_l, \mu_r)$ is an $H$-superalgebra. We note that if $H$ is trivial,
i.e. $H = \C$, then an $H$-superalgebra is merely a superalgebra.

Next we define the algebra $D_G$ precisely the same way that we defined it in
\S\ref{QuantumGroupoids}. Once again, there is a natural embedding of
$M_G$ in $D_G$. We can then give the definition for an
\emph{$H$-superbialgebra}. 
As expected, for this, we only need a \emph{coproduct}, i.e. a
homomorphism of $H$-superalgebras ${\Delta : A \rightarrow A \otimes A}$ and a
\emph{counit}, i.e. a homomorphism of $H$-superalgebras ${\epsilon : A
\rightarrow D_G}$. We require that $\Delta$ be coassociative, i.e. ${(\Delta
\otimes id_A) \circ \Delta = (id_A \otimes \Delta) \circ \Delta}$, and that the
counit $\epsilon$ satisfy the \emph{counit axiom}: 
\[ (\epsilon \otimes id_A) \circ \Delta = (id_A \otimes \epsilon) \circ \Delta
= id_A. \]
\noindent
In short, we say that $\mathcal{A}_H = (A, H, \mu_l, \mu_r ,\Delta, \epsilon)$
is an \emph{$H$-superbialgebra} if $(A, H, \mu_l, \mu_r)$ is an
$H$-superalgebra, and $\Delta$ and $\epsilon$ satisfy the required conditions
above. We note that if $H$ is trivial, i.e. $H = \C$, then an
$H$-superbialgebra is merely a superbialgebra.

To obtain an $H$-Hopf superalgebra, we only need to add to the above an
\emph{$H$-antipode}, i.e. an antiautomorphism ${S : A \rightarrow A}$
satisfying
the technical condition:
\begin{eqnarray*}
a_{(1)}S(a_{(2)}) &=& \mu_l(\epsilon(a)1), \\
S(a_{(1)})a_{(2)} &=& \mu_r(\sigma_{\alpha}\epsilon(a)1),
\end{eqnarray*}
\noindent
where $a \in A_{\alpha,\beta}$, $\epsilon(a)1$ is the result of applying the
difference
operator $\epsilon(a) \in D_G$ to $1 \in M_G$ and we used the Sweedler notation
associated to $\Delta(a)$. With this last ingredient, we can finally define an
\emph{$H$-Hopf superalgebra}: It is an $H$-superbialgebra with an $H$-antipode. 
Once again in the case when $H = \C$, an $H$-Hopf superalgebra is only a
super Hopf algebra.\footnote{Here, we may mention the idea of
\emph{bosonization}, which was introduced by Radford in \cite{Rad}. For
another exposition, the reader may look at \cite{Maj}. The connection between
Hopf and super Hopf structures may be viewed as a special case of
bosonization; see \cite[Sec. 7]{Gel}. The main idea of bosonization with
respect to super Hopf algebras is as follows: For many practical
purposes we can just as well work with regular Hopf algebras instead of super
Hopf algebras. More precisely, the tensor category of representations of an
ordinary Hopf algebra $H$ with a grouplike element $g$ with $g^2 = 1$ is
equivalent to the tensor category of representations of a related super Hopf
algebra $H_s$ with a grouplike odd element $g_s$ with $g_s^2 = 1$, and the
correspondence between such pairs $(H,g)$ and $(H_s,g_s)$ is one-to-one. The
significance and the implications of this categorical equivalence for our
context need further investigation.}

Now, we are ready to move on to the superization of the constructions of \S\ref{CategoricalSubsection}. We first define super dynamical representations: A \emph{super
dynamical representation} of an $H$-superalgebra $A$ should be a pair
$(W,\pi_W)$, where $W$ is a diagonalizable $H$-module and $\pi_W : A
\rightarrow D_{G,W}$ is a homomorphism of $H$-superalgebras. Here, as before,
$H$ is taken to be a commutative, cocommutative, finitely generated Hopf
algebra over $\C$. $G$, $M_G$, $D_G$, $D_{G,W}$ and $D^{\alpha}_{G,W}$ are
defined as above. Note that, this time, $W$ lives in the category of super vector
spaces.

Next, we will follow \S\ref{CategoricalSubsection} and superize the construction of tensor products of dynamical representations in order to obtain the right kind of category. In particular, we will let the tensor product $V \hat{\otimes} W$ of two dynamical representations $V$ and $W$ of an $H$-superbialgebra $A$ to be the usual tensor product $V \otimes W$ in the category of super vector spaces, and we will define the map $\pi_{V
\hat{\otimes} W} : A \rightarrow D_{G,V \hat{\otimes} W}$ by setting $\pi_{V
\hat{\otimes} W}(a) = \theta_{V,W} \circ (\pi_V \otimes \pi_w) \circ
\Delta(a)$. Here $\theta_{V,W}$ is the natural embedding, as a
$H$-superalgebra, of $D_{G,V} \hat{\otimes} D_{G,W}$ into $D_{G, V
\hat{\otimes} W}$ given by:
\[ f\sigma_{\beta} \otimes g \sigma_{\delta} \longmapsto (f \hat{\otimes} g)
\sigma_{\delta} \]
\noindent
where $\sigma_{\beta}$ is the automorphism of $M_G$ (embedded naturally inside
$D_G$), given by $\sigma_{\beta}(f) (\lambda) = f(\lambda + \beta)$ and
$\sigma_{\delta}$ is defined in a similar manner. The function $f
\hat{\otimes} g$ is determined by: 
\[f \hat{\otimes} g (\lambda)(v \otimes w) = f(\lambda-\mu)(v) \otimes
g(\lambda)(w) \]
\noindent
if $g(\lambda)(w)$ has weight $\mu$. This construction will give us a tensor
category.\footnote{If we want rigidity for our category, we need to construct left and right duals, and this can be done, with the additional hypothesis that $A$ be an $H$-Hopf superalgebra. We will once again choose not to go into the duality problem.} 

Let $V_d$ be the tensor category of super dynamical representations of an
$H$-superbialgebra or an $H$-Hopf superalgebra, and assume that we are given a
tensor functor ${\mathcal{F} : \mathcal{B} \rightarrow V_d}$. Then using any
object $X$ of $\mathcal{B}$, we can construct a $V_d$-automorphism $R_s =
{R_s(\mathcal{B},\mathcal{F},X)}$ of $\mathcal{F}(X) \otimes \mathcal{F}(X)$
by setting $R_s =T_s F (\beta_{X,X})$, where ${\beta_{X,Y} : X
\otimes_{\mathcal{B}} Y \rightarrow Y \otimes_{\mathcal{B}} X}$ is the
braiding of the category $\mathcal{B}$ and $T_s$ is the usual super twist in
the category of super vector spaces. It is a simple exercise to show that this
$R_s$ satisfies the quantum dynamical Yang-Baxter equation (QDYBE) and hence
is a \emph{super dynamical $R$-matrix}. Conversely it can be shown that given
a super dynamical $R$-matrix $R_s$, i.e. a solution $R_s$ to the QDYBE, we can
find appropriate categories $\mathcal{B}$ and $V_d$, a tensor functor
${\mathcal{F} : \mathcal{B} \rightarrow V_d}$ and a particular object $X$ of
$\mathcal{B}$ such that $R_s$ is the $V_d$-automorphism associated to the
triple ${(\mathcal{B},\mathcal{F},X)}$ in the same manner as in \S\ref{CategoricalSubsection}; this is a
straight-forward extension of Theorem 3.1 of \cite{EV2}.

So far, we see that almost all our constructions are straight-froward superizations of those from the non-graded theory. 

Finally we are in a position to define \emph{super dynamical quantum groups}:
These are going to be those $H$-Hopf superalgebras which can be obtained from
super dynamical $R$-matrices satisfying the \emph{super Hecke condition}. We
will discuss this condition and the associated construction in more detail in
the next subsection. 

We end this subsection with a note about our choice of terminology: In this 
paper we preferred to use the terms \emph{$H$-superbialgebra}, 
\emph{$H$-Hopf superalgebra}, and \emph{super dynamical quantum group} 
instead of the perfectly acceptable alternatives \emph{$H$-superbialgebroid}, 
\emph{$H$-Hopf superalgebroid}, and \emph{super} \emph{dynamical} 
\emph{quantum} \emph{groupoid} in order to minimize the number of terms which 
end with $-oid$. A plausible argument for the second set of terms could be that
these signify much more clearly that they are special types of 
\emph{superbialgebroids}, \emph{Hopf superalgebroids} or \emph{super quantum
groupoids}, but since we did not actually define these latter terms, we feel
comfortable with our choices.

\subsection{Super Dynamical $R$-matrices}
\label{ResultsSuperGroups}

Here we present some constructions related to super dynamical quantum groups.
We start with a discussion of the QDYBE and its solutions in the super setting. We follow closely \cite{ES2, EV2}.

We restrict ourselves to the following setup: Let $\h$ be the standard
Cartan subsuperalgebra of $gl(m,n)$. This is the set of all diagonal matrices
from $gl(m,n)$. Let $V$ be a semisimple $\h$-module. Recall that the quantum
dynamical Yang-Baxter equation is:
\begin{equation*}
\tag*{(\ref{QDYBE})}
R^{12}(\lambda- h^{(3)}) R^{13}(\lambda) R^{23}(\lambda-h^{(1)}) =
R^{23}(\lambda) R^{13}(\lambda - h^{(2)}) R^{12}(\lambda).
\end{equation*}
\noindent 
We will say that a solution $R : \h^* \rightarrow End_{\h}(V \otimes V)$ to the
QDYBE is a \emph{super dynamical $R$-matrix} if it is an invertible
meromorphic function. For simplicity, we will only consider even solutions.

Rescaling the QDYBE by $\lambda \mapsto \tfrac{\lambda}{\gamma}$ we get the
\emph{quantum dynamical Yang-Baxter equation with step} $\gamma$:
\begin{equation*}
R^{12}(\lambda- \gamma h^{(3)}) R^{13}(\lambda) R^{23}(\lambda-\gamma h^{(1)})
=
R^{23}(\lambda) R^{13}(\lambda - \gamma h^{(2)}) R^{12}(\lambda).
\end{equation*}
\noindent 
Invertible solutions of the QDYBE with step $\gamma$ will also be called
\emph{super dynamical $R$-matrices}. 

If $R : \h^* \rightarrow End_{\h}(V \otimes V)[[\gamma]]$ is a continuous
family of $\h$-invariant meromorphic functions of the form $R = 1 - \gamma r +
O(\gamma^2)$ satisfying the QDYBE with step $\gamma$, then a simple computation
will show that $r$ is an $\h$-invariant super dynamical $r$-matrix satisfying
the CDYBE. This follows from a straight-forward superization of Proposition
3.1 of \cite{ES2}. We naturally
call $r$ the \emph{semiclassical limit} of $R$ and $R$ a \emph{quantization} of
$r$.

Now we consider the case when $V = V_{\overline{0}} \oplus V_{\overline{1}}$ is
the standard vector representation of $\h$. In other words, $V$ is a super
vector space with even dimension $m$ (i.e. $\dim_{\C}(V_{\overline{0}}) =
m$) and odd dimension $n$ (i.e. $\dim_{\C}(V_{\overline{1}}) = n$). 

Let
$\{h_{\overline{0},1}, \cdots, h_{\overline{0},m}, h_{\overline{1},1}, \cdots,
h_{\overline{1},n}\}$ be the standard basis for $\h$ and let
${\lambda}_{\overline{0},1}, \cdots, {\lambda}_{\overline{0},m},
{\lambda}_{\overline{1},1}, \cdots, {\lambda}_{\overline{1},n}$ be the
corresponding coordinate functions on $\h^*$. Note that each 
$h_{\overline{s},i}$ will be even as all of $\h$ lies inside
$gl(m,n)_{\overline{0}}$. Let $V_{\overline{s},i}$, $s =
0,1$, be the one-dimensional weight subspaces of $V$ of weight
$\omega_{\overline{s},i}$ such that $(\omega_{\overline{s},i},
h_{\overline{t},j}) = \delta_{s,t}\delta_{i,j}$. 
%Using this notation, we can
%see that 
Then, the tensor product $V \otimes V$ will have weight subspaces of the
form $V_{\overline{s},i} \otimes V_{\overline{s},i}$ and $(V_{\overline{s},i}
\otimes V_{\overline{t},j}) \oplus (V_{\overline{t},j} \otimes
V_{\overline{s},i})$.

We will say that a super dynamical $R$-matrix $R$ satisfies the
\emph{super Hecke condition}\footnote{Analogous to the non-graded case, the super Hecke condition is a quantum version of Equation
\eqref{superdynamicalunitarity}, the generalized unitarity condition. In
particular, if a continuous family of super dynamical $R$-matrices of the form
$R = 1 - \gamma r + O(\gamma^2)$ satisfies the super Hecke condition with
parameter $q =1$, then the semiclassical limit $r$ satisfies the
unitarity condition $r + T_s(r) = 0$.}, with parameter $q \in \C^*$ if the
eigenvalues of $T_sR$ are $1$ on $V_{\overline{s},i} \otimes
V_{\overline{s},i}$ and $1, -(-1)^{st}q$ on $(V_{\overline{s},i} \otimes
V_{\overline{t},j}) \oplus (V_{\overline{t},j} \otimes V_{\overline{s},i})$.
Note that here $T_s$ stands for the matrix form of the super twist, but since
$R$ is assumed to be even, it may be replaced simply by the permutation
matrix.

We will say that a super dynamical $R$-matrix has \emph{zero weight} if we
have:
\begin{equation}
\label{zeroweightsuper} 
[R^{ij}(\lambda), h \otimes 1 + 1 \otimes h] = 0 
\end{equation}
\noindent
for all $i,j = 1,2,3$, $i \neq j$ and $h \in \h$, $\lambda \in \h^*$. As in the
case of super dynamical $r$-matrices in Subsection
\ref{GeneralizingSchiffmann}, the \emph{zero weight condition}, i.e. Equation 
\eqref{zeroweightsuper}, is equivalent to
$\h$-invariance. 

We can
easily see that a zero-weight super dynamical $R$-matrix $R(\lambda)$
satisfying the super Hecke condition with parameter $q$ has to be of the form:
\[ R(\lambda) = \sum_a E_{aa} \otimes E_{aa} + \sum_{a \neq b}
\alpha_{ab}(\lambda) E_{aa} \otimes E_{bb} + \sum_{a \neq b}
\beta_{ab}(\lambda)E_{ab} \otimes E_{ba},
\]
\noindent
where $E_{ab}$ is the elementary matrix coresponding to the $(a,b)$th matrix
entry, and $\alpha_{ab}$ and $\beta_{ab}$ are certain meromorphic functions
$\h^* \rightarrow \C$.

In \S\ref{QuantumGroupoids}, we noted in passing that 
dynamical $R$-matrices satisfying the (non-graded) Hecke condition can be
used to define $H$-Hopf algebroids. Now we will explicitly superize the
relevant constructions from \cite{ES2, EV2} in order to see how we can obtain
an $H$-Hopf superalgebra from a super dynamical $R$-matrix satisfying the
super Hecke condition. 

For this, we start with our commutative Lie superalgebra $\h$, the standard
Cartan subsuperalgebra of $gl(m,n)$. Recall that this happens to be the set of
all diagonal matrices from $gl(m,n)$. Once again, we let $V = V_{\overline{0}}
\oplus V_{\overline{1}}$ be a semisimple $\h$-module (living in the category
of super vector spaces). 

Set $\dim_{\C}(V_{\overline{0}}) = M$ and
$\dim_{\C}(V_{\overline{1}}) = N$. Pick an ordered homogeneous basis for $V$:
$\{v_1,\cdot,v_{M+N}\} = \{v_{\overline{0},1}, \cdots, v_{\overline{0},M}, 
v_{\overline{1},1}, \cdots, v_{\overline{1},N}\}$ with $v_{\overline{p},i} \in
V_{\overline{p}}$, $p= 0,1$.

Now let $R : \h^* \rightarrow End_{\h}(V \otimes V)$ be a meromorphic function
satisfying the zero weight condition such that $R(\lambda)$ is invertible for
generic $\lambda$. We denote by $H$ the super Hopf algebra associated to $\h$;
in other words, $H$ is the commutative cocommutative super Hopf algebra which
is the universal enveloping superalgebra of $\h$. Clearly $V$ is also an
$H$-module. Let $G$, $M_G$ and $D_G$ be defined as before. In other words, $G =
Spec(H)$, $M_G$ is the field of meromorphic functions on $G$, and $D_G$ is the
algebra of all difference operators on $M_G$.

We first define an $H$-superbialgebra $\mathcal{A}_H(R) = (A_R, H, \mu_l^R,
\mu_r^R ,\Delta_R, \epsilon_R)$ associated to $R$. For this, we need an
$H$-superalgebra $(A_R, H, \mu_l^R, \mu_r^R)$ and two maps $\Delta_R$ and
$\epsilon_R$ satisfying the required conditions from Subsection
\ref{DefinitionsSuperGroups}. To describe $A_R$, we start with the
superalgebra generated freely by $M_G \otimes M_G$. We follow \cite{EV2} in
their notations and denote the first copy of $M_G$ by the superscript $^{(1)}$
and the second by $^{(2)}$: If $f \in M_G$ is any function, then we denote the
image of $f$ in the first (respectively, the second) copy of $M_G$ by $f(\lambda^{(1)})$ (respectively, by $f(\lambda^{(2)})$).

Next we put in some additional generators $L_{i,j}$, $L^{-1}_{i,j}$, for $i,j$
between $1$ and $M+N$. Intuitively, each $L_{i,j}$ corresponds to the $(i,j)$th
matrix entry of a particular operator: 
\[ L = \sum_{a,b = 1}^{M+N} E_{ab} \otimes L_{a,b}, \] 
\noindent
and each $L^{-1}_{i,j}$ corresponds to the $(i,j)$th matrix entry of the
inverse operator $L^{-1}$. The parity of these generators is determined once a
homogeneous basis for $V$ is picked. Hence, $L_{i,j}$ and $L^{-1}_{i,j}$ will
be even if $(1\le i,j \le M)$ or $(M+1 \le i,j \le M+N)$, and they will be odd
otherwise. 

Now we will define $A_R$ to be the quotient of this superalgebra by several
relations. First we will encode the invertibility of $L$ by requiring that:
\[ L L^{-1} = L^{-1} L = 1 \]
\noindent
and enforce the QDYBE by requiring:
\begin{equation}
\label{QDYBEwithL}
R^{12}(\lambda^{(1)})L^{13}L^{23} = \quad :L^{23} L^{13}
R^{12}(\lambda^{(2)}): 
\end{equation}
\noindent
Note that in the equation above, we are using the normal order notation ``::".
In
other words, if we define the matrix entries of $R(\lambda)$ by:
\[ R(\lambda)(v_{a}\otimes v_{b}) = \sum_{c,d=1}^{M+N}
R^{ab}_{cd}(\lambda)v_c\otimes v_d, \]
\noindent
then the above version of the QDYBE may be rewritten as follows:
\[ \sum_{x,y=1}^{M+N} R^{xy}_{ac}(\lambda^{(1)})L_{x,b}L_{y,d} =
\sum_{x,y=1}^{M+N}
R^{bd}_{xy}(\lambda^{(2)})L_{c,y}L_{a,x}. \]

We also require the following relations:
\begin{eqnarray*}
f(\lambda^{(1)})L_{\alpha,\beta} &=& L_{\alpha,\beta} f(\lambda^{(1)} + \alpha)
\\
f(\lambda^{(2)})L_{\alpha,\beta} &=& L_{\alpha,\beta} f(\lambda^{(2)} + \beta)
\end{eqnarray*}
\noindent
along with
\[ [f(\lambda^{(1)}),g(\lambda^{(2)})] =0.\]
\noindent
Here, $L_{\alpha,\beta}$, $\alpha, \beta \in G$, are the weight
components of $L$ with respect to the $G$-bigrading on $End(V)$ corresponding
to the $G$-weight decomposition $V = \oplus_{\alpha \in G}
V_{\alpha}$.\footnote{Equivalently we can extract the weight decomposition of $L$ from
the natural $\h$-bigrading on $End(V)$
corresponding to the weight decomposition $V = \oplus_{\alpha \in \h^*}
V_{\alpha}$ if we prefer to think only in terms of the Lie superalgebra $\h$
and its associates.}

The last few relations make more sense once we define the moment maps $\mu_l^R$
and $\mu_r^R : M_{G} \rightarrow (A_R)_{0,0}$:
\[ \mu_l^R(f) = f(\lambda^{(1)}) \qquad \textmd{and}\qquad \mu_r^R(f) =
f(\lambda^{(2)}). \]
\noindent
%Here, the subscript on $(A_R)_{0,0}$ denotes the $G$-bigrading. 
We also require that $\mu_l^R$ and $\mu_r^R$ will preserve the
$\Z/{2\Z}$-grading. Now we can
see that the two relations about $f(\lambda^{(i)})L_{\alpha,\beta}$ are merely
the defining relations of the moment maps.

At this point, it is clear that the quadruple $(A_R, H, \mu_r^R, \mu_l^R)$ is an
$H$-superalgebra. To obtain the $H$-superbialgebra that we want, we only need
to define a suitable coproduct ${\Delta_R : A_R \rightarrow A_R \otimes A_R}$
and an appropriate counit ${\epsilon_R : A_R \rightarrow D_G}$. The
definitions from \cite{ES2, EV2} will work for us: 
\[ \Delta_R(L) = L^{12}L^{13}, \quad 
\Delta_R(L^{-1}) = (L^{-1})^{13}(L^{-1})^{12}, \]
\noindent
and:
\[ \epsilon_R(L_{\alpha,\beta}) = \delta_{\alpha,\beta} id_{V_{\alpha}} \otimes
\sigma^{-1}_{\alpha}, \quad \epsilon_R((L^{-1})_{\alpha,\beta}) =
\delta_{\alpha,\beta} id_{V_{\alpha}} \otimes \sigma^{-1}_{\alpha}.\]
\noindent
The arguments in the proofs of Propositions 4.2 and 4.3 from \cite{EV2} can
then be superized in a straight-forward manner to show that the defining
relations of $A_R$ are invariant under $\Delta_R$ and are annihilated by
$\epsilon_R$, that the counit satisfies the counit axiom, and finally that
$\mathcal{A}_H(R) = (A_R, H, \mu_l^R, \mu_r^R ,\Delta_R, \epsilon_R)$ is an
$H$-superbialgebra.

Now to get an $H$-Hopf superalgebra, we need  
an appropriate $H$-antipode.
Once again we follow the arguments of \cite{EV2} and see that they are easily
generalized to the super case. 
The proofs of Propositions 4.4, 
4.5 and 4.6 from \cite{EV2} carry over almost directly. We can show that, 
if we start with a continuous family of super dynamical $R$-matrices 
$R_{\gamma}(\lambda)$ with step $\gamma$ such that $R_0 = 1$, then there is 
a unique map $S_R$ on $A_R$ which sends $L$ to $L^{-1}$ and makes $A_R$ into 
an $H$-Hopf superalgebra $\mathcal{A}_H(R) = (A_R, H, \mu_l^R, \mu_r^R,
\Delta_R, \epsilon_R, S_R)$.

In the above constructions, we did not use the QDYBE. In fact 
dropping the QDYBE, or more accurately Equation \eqref{QDYBEwithL}, would 
still give us an $H$-superbialgebra. However, we want more than just an 
$H$-superbialgebra. In \cite{ES2}, Etingof and Schiffmann show 
that the QDYBE is a sufficient condition to ensure the existence of 
at least one dynamical representation of $\mathcal{A}_H(R)$. Their 
argument carries over to the super case. If $R$ is a super 
dynamical $R$-matrix, then we can obtain a super dynamical representation 
$(V,\pi_V)$ of $\mathcal{A}_H(R)$ on $V$ by
first defining a map 
${\pi_V^0 : A_R \rightarrow Hom(V, V \otimes M_G)}$ via
${\pi_V^0(\lambda)} = R(\lambda)$, and, then, taking $\pi_V$ as 
the unique extension of $\pi_V^0$. 

We end this section with a short remark about the super Hecke condition. 
Notice that we also did not make use of the super Hecke condition in the 
construction above. In the non-graded case, when we restrict the types of 
algebras that can be constructed by the above technique, the Hecke condition 
comes up naturally. In fact, it is shown explicitly in \cite{ES2} that, if 
the (non-graded) dynamical $R$-matrix satisfies the (non-graded) Hecke 
condition, then the associated algebra $A_R$ may be viewed as a deformation 
of a function algebra of a matrix group. We leave the discussion of the 
superization of this result to the next section.

\section{Open Questions	}
\label{SectionConclusion}

In this paper, we aimed to discuss both what has been accomplished so far in the
theory of super quantum groups and the path still lying ahead. Our approach
definitely favored Hopf algebroids as the main framework, and the emphasis was 
clearly on the dynamical picture. 

We now end this paper with a brief discussion of the results presented and the
questions that still need to be answered. 

As we have seen, many constructions related to dynamical quantum groups can
naturally be extended to the super case. Many definitions are the expected
super analogues of the non-graded ones. We have presented the super versions of 
the major classification results in the classical picture, proving explicitly 
that the constructions still make sense. For the quantum picture, we have only 
just begun: we made the necessary definitions and extended only the most basic 
of the constructions. 

Below we list some remarks and questions that we think may be relevant for the 
further development of the theory of super quantum groups. This list may also 
be viewed as our plan of action for the near future:

\textbf{(1).} In \S\ref{ResultsSuperGroups}, we proposed a superization of the
Hecke 
condition that we think best generalizes the non-graded version by 
considering it as the right way to quantize Equation 
\eqref{superdynamicalunitarity}, the generalized unitarity condition. A good test 
to see whether we made the correct definition would be to check if the  
superalgebra $A_R$ associated to a super dynamical $R$-matrix satisfying the 
super Hecke condition may be viewed as a deformation superalgebra of the 
appropriate type. In fact, the relevant results from \cite{EV2} (Theorems 6.1,
6.2, 
6.3 and Proposition 6.1) seem, at a first glance, amenable to straight-forward 
superization, and we expect that this should not be hard to verify.

\textbf{(2).} In \cite{EV2}, Etingof and Varchenko classify all (non-graded)
dynamical
$R$-matrices that satisfy the Hecke condition. If the super Hecke condition is correct, then a natural next step would be to extend 	this result to the super case. 

\textbf{(3).} There is a nice path to a partial solution of the quantization problem in the super setting. In \cite{EV2}, all (non-graded) zero-weight dynamical $r$-matrices (with any coupling constant) are 
explicitly quantized. We have seen in \S\ref{SuperZeroWeight} that the 
classification paradigm for the zero-weight dynamical $r$-matrices in the 
non-graded case can be superized in a natural manner. Therefore, one could
expect
that the solution of the quantization problem for the zero-weight super 
dynamical $r$-matrices will be a natural superization of the solution of the 
non-graded problem.

\textbf{(4).} The general problem of quantization remains open. The extension
to the 
graded world of the main quantization result from \cite{ESS}, the general 
constructive quantization of all classical dynamical $r$-matrices which fit
Schiffmann's classification (\S\ref{SchiffmannResult}) is very 
important.

\textbf{(5).} The complications mentioned in \S\ref{SuperSolutionsNonDynamical} which made the Belavin-Drinfeld classification results hard to superize are still there in the dynamical case. These still need to be addressed.
Hence, the classification problem is still open as well.

\end{document}